\documentclass[11pt,a4paper]{article}

\usepackage[latin1]{inputenc}
\usepackage{tikz}
\usetikzlibrary{trees}

\usepackage[tbtags]{amsmath}    
\usepackage{amsfonts,amssymb,amsthm,euscript,makeidx,pifont,boxedminipage,multicol,fullpage}
\usepackage{color}
\usepackage{yfonts}
\usepackage[authoryear]{natbib}

\def\be{\begin{eqnarray}}
\def\ee{\end{eqnarray}}

\def\b*{\begin{eqnarray*}}
\def\e*{\end{eqnarray*}}

\newtheorem{Theorem}{Theorem}[section]
\newtheorem{Lemma}[Theorem]{Lemma}
\newtheorem{Proposition}[Theorem]{Proposition}
\newtheorem{Corollary}[Theorem]{Corollary}
\newtheorem{Definition}[Theorem]{Definition}
\theoremstyle{definition}
\newtheorem{Remark}[Theorem]{Remark}
\newtheorem{Example}[Theorem]{Example}

\newcommand{\rmi}{{\rm (a)$\>\>$}}
\newcommand{\rmii}{{\rm (b)$\>\>$}}
\newcommand{\rmiii}{{\rm (c)$\>\>$}}
\newcommand{\rmiv}{{\rm (d)$\>\>$}}

\def \Rbrack{]\!]}
\def \Lbrack{[\![}

\def \E{\mathbb{E}}
\def \F{\mathbb{F}}
\def \G{\mathbb{G}}
\def \H{\mathbb{H}}

\def \P{\mathbb{P}}
\def \Q{\mathbb{Q}}
\def \R{\mathbb{R}}

\def \ff{\mathbb{F}}

\def\Ac{{\cal A}}
\def\Bc{{\cal B}}
\def\Cc{{\cal C}}

\def\Ec{{\cal E}}
\def\Fc{{\cal F}}
\def\Gc{{\cal G}}
\def\Hc{{\cal H}}
\def\hcal{{\mathfrak h}}

\def\Lc{{\cal L}}
\def\Mc{{\cal M}}

\def\Pc{{\cal P}}

\def\Rc{{\cal R}}

\def\Tc{{\cal T}}

\def\Xc{{\cal X}}
\def\Yc{{\cal Y}}

\def \id{1\!\!1}
\def \Law{\mathcal{L}}

\def \Om{\Omega}
\def \om{\omega}
\def \Omb{\overline{\Om}}
\def \omb{\overline{\om}}
\def \omh{\widehat{\om}}
\def \y{y}
\def \Y{Y}
\def \fu{f}
\def \Yc{\mathcal Y}
\def \Sh{\widehat S}

\def \pib{\overline{\pi}}
\def \pih{\widehat{\pi}}
\def \pit{\overline{\pih}}

\def \xib{\overline{\xi}}

\def \eps{\varepsilon}

\def \0{\mathbf{0}}
\def \1{\mathbf{1}}
\def \x{\times}
\def \no{\lambda}
\def \No{\Lambda}
\def \Noh{\widehat{\No}}
\def \noh{\widehat{\no}}

\def \Fcb{\overline{{\cal F}}}
\def \Gcb{\overline{{\cal G}}}
\def \Fbb{\overline{\F}}
\def \Gbb{\overline{\G}}

\def \Mct{\overline{\Mch}}
\def \Ect{\overline{\widehat{{\cal E}}}}
\def \Ectt{\widetilde{{\cal E}}}

\def \Qt{\overline{\Qh}}

\def \Pcb{\overline{\Pc}}

\def \Pb{\overline{\P}}

\def \Hcb{\overline{\Hc}}
\def \Hb{\overline{H}}
\def \Hh{\widehat{H}}

\def \is{\circ}

\def\esup{\mathop{\rm ess \, sup}}
\def\Bf{\mathfrak{P}}
\def\Rb{\overline{\R}}
\def\T{\mathbb{T}}
\def\Qb{\overline{\Q}}
\def\Mcb{\overline{\Mc}}
\def\Mce{\underline{\Mc}}
\def\Ecb{\overline{\Ec}}
\def\Ech{\widehat{\Ec}}
\def\Lcb{\overline{\Lc}}
\def\Ccb{\overline{\Cc}}
\def\NA{\mathrm{NA}}

\def\Omh{\widehat{\Om}}
\def\Xh{\widehat{X}}
\def\Pch{\widehat{\Pc}}
\def\Mch{\widehat{\Mc}}

\def\Pbh{\widehat{\mathbb P}}
\def\Qh{\widehat{\Q}}
\def\Fh{\widehat{\F}}
\def\Gh{\widehat{\G}}
\def\Fch{\widehat{\Fc}}
\def\Gch{\widehat{\Gc}}

\def\xih{\widehat{\xi}}

\def\Psih{\widehat{\Psi}}

\def\Hh{\widehat{H}}
\def\Hch{\widehat{\mathcal H}}
\def\omh{\widehat{\om}}

\def\Sh{\widehat{S}}

\def\fcl{\Upsilon}
\def\fclb{\overline\fcl}
\def\fclh{\widehat\fcl}

\title{
Robust pricing--hedging duality for American options in discrete time financial markets
\footnote{We are grateful to Bruno Bouchard, David Hobson, Zhaoxu Hou and Monique Jeanblanc for useful comments and suggestions.
SD and XT gratefully acknowledge the financial support of the ERC 321111 Rofirm, the ANR Isotace, and the Chairs Financial Risks (Risk Foundation, sponsored by Soci\'et\'e G\'en\'erale) and Finance and Sustainable Development (IEF sponsored by EDF and CA). AA and JO thankfully acknowledge support from the ERC under the EU's FP7/2007-2013 (ERC grant no.\ 335421). JO is grateful to St John's College for its support.}}

\author{
	Anna Aksamit\thanks{University of Oxford aksamit@maths.ox.ac.uk and obloj@maths.ox.ac.uk}
	\and Shuoqing Deng\thanks{University of Paris-Dauphine, PSL Research University, CNRS, UMR [7534], CEREMADE. deng@ceremade.dauphine.fr and 
	tan@ceremade.dauphine.fr}
	\and Jan Ob\l\'{o}j\footnotemark[2]
	\and Xiaolu Tan\footnotemark[3]
}

\date{\today}

\begin{document}

\maketitle

\abstract{
We investigate pricing--hedging duality for American options in discrete time financial models where some assets are traded dynamically and others, e.g.\ a family of European options, only statically.
In the first part of the paper we consider an abstract setting, which includes the classical case with a fixed reference probability measure as well as the robust framework with a non-dominated family of probability measures. Our first insight is that by considering a (universal) enlargement of the space, we can see American options as European options and recover the pricing--hedging duality, which may fail in the original formulation. This may be seen as a weak formulation of the original problem. 
Our second insight is that lack of duality is caused by the lack of dynamic consistency and hence a different enlargement with dynamic consistency is sufficient to recover duality: it is enough to consider (fictitious) extensions of the market in which all the assets are traded dynamically. 
In the second part of the paper we study two important examples of robust framework: the setup of \cite{BN13} and the martingale optimal transport setup of \cite{BHLP}, and show that our general results apply in both cases and allow us to obtain pricing--hedging duality for American options.
}

\vspace{2mm}


\vspace{2mm}


\section*{Introduction}
Robust approach to pricing and hedging has been an active field of research in mathematical finance over the recent years. In this approach, instead of choosing one model, one considers superhedging simultaneously under a family of models, or pathwise on a set of feasible trajectories.  
It generalizes the classical approach where one only considers models which are absolutely continuous with respect to a fixed reference probability measures $\P$.
In such setting, absence of arbitrage is known to be equivalent to the existence of an equivalent martingale measure, result known as the first fundamental theorem of asset pricing, see e.g. \cite{DeSch}, \cite{FS04}. 
	When the market is complete --- i.e.\ when every contingent claim can be perfectly replicated by a self--financing trading strategy ---
	the equivalent martingale measure $\Q$ is unique and option prices are given by their replication cost, which is equal 
	to the expected value of the discounted payoff under $\Q$. 
	In an incomplete market, when a perfect replication strategy does not exist,
	a safe way of pricing is to use the minimum super--replication cost of the option.
	Using duality techniques, this minimal super--replication price can be related to the pricing problem and expressed as the supremum of expectations of the discounted payoff over all martingale measures equivalent to $\P$. 
		
One of the challenges in a robust approach is to extend this dual relationship to non-dominated (robust) context.
In continuous time models under \emph{volatility uncertainty}, such pricing--hedging duality results have been obtained by, among many others, \citet{DM06}, \cite{STZ}, \cite{NN12}, \cite{PRT}. In discrete time, general pricing--hedging duality was shown in e.g. \cite{BN13} and \cite{BFM15}. 
Importantly, in a robust setting one often wants to include further market instruments which may be available for trading. 
In a setup which goes back to the seminal work of  \cite{Hlookback}, one often considers dynamic trading in the underlying and  static trading, i.e.\ buy and hold strategies at time zero, in some European options, often call or put options with a fixed maturity. Naturally, such additional assets constraint the set of martingales measures we may use for pricing. General pricing--hedging duality results in variants of this setting, both in continuous and in discrete time, can be found in e.g.\ \cite{ABPS}, \cite{BFHMO16}, \cite{BHLP}, \cite{DS14a}, \cite{HouObloj}, \cite{GTT}, \cite{TT} and we refer to the survey papers \cite{Hreview, Obloj} for more details.

The main focus in the literature so far has been on the duality for (possibly exotic) European payoffs. However, more recently, some focus was put on American options.  \cite{CoxHoeggerl} 
	studied the necessary (and sufficient in some cases) conditions on the American put option prices for the absence of arbitrage.
\cite{Dolinsky14} studied game options (including American options) in a non--dominated discrete--time market, but without statically traded options allowed for super--replication. \cite{Neuberger} considered a discrete--time, discrete space market with presence of statically traded European vanilla options. He observed that the superheding price for an American option may be strictly larger than the supremum of the expected (discounted) payoff over all stopping times and all (relevant) martingale measures. We refer to this situation as \emph{duality gap}. In \cite{Neuberger}, the pricing--hedging duality is then restored by using a \emph{weak} dual formulation. This approach was further exploited, with more general results, in \cite{HobsonNeuberger}. \cite{Bay} studied the same superhedging problem in the setup of \cite{BN13} but only considered \emph{strong} stopping times in their dual formulation, which leads to a duality gap in general.
More recently, and in parallel to an earlier version of this paper, \cite{BayZhou} proved a duality result by considering \emph{randomized} models, under some regularity and integrability conditions on the payoff functions.

Motivated by the above works, we endeavour here to understand the fundamental reasons why pricing--hedging duality for American options holds or fails and offer systematic reasons to mend it in the latter case. 
We derive two main general results which we then apply to various specific contexts, both classical and robust.
Our first insight is that by considering a (universal) enlargement of the space, namely the time--space product structure, we can see an American option as a European option and recover the pricing--hedging duality, which may fail in the original formulation. 
This may be seen as a weak formulation of the dual (pricing) problem and leads to considering a large family of stopping times. 
	This formulation of the dual problem is similar in spirit to \cite{Neuberger, HobsonNeuberger, BayZhou},
	but our approach leads to a duality results in a more general setting, and/or under more general conditions,  see Remark \ref{com:hob} and Subsection \ref{com:bay} and also \cite{HobsonNeuberger2}.
	Our second main insight is that the duality gap is caused by the failure of dynamic programming principle. 
	To recover the duality, under the formulation with \emph{strong} stopping times, it is necessary and sufficient to consider an enlargement which restores the dynamic consistency: it is enough to consider (fictitious) extensions of the market in which all the assets are traded dynamically.
	As a byproduct, we find that the dynamic trading strategy on options and the classical semi-static strategy lead to the same superhedging cost in various settings.

The first part of paper, Section \ref{duality}, presents the above two main insights in a very general discrete time framework which covers both classical (dominated) and robust (non--dominated) settings. In the second part of the paper, we apply our general results in the context of two important examples of the robust framework: the setup of \cite{BN13} in Section \ref{sec:main} and the martingale optimal transport setup of \cite{BHLP} in Section \ref{sec:mot}. We obtain suitable pricing--hedging duality for American options in both setups. 
In the latter case of martingale optimal transport, there is an infinity of assets to consider and we use measure valued martingales to elegantly describe this setting.

\begin{Example}
\label{ex:intro}
We conclude this introduction with a motivating example showing that the pricing--hedging duality may fail in presence of static trading instruments and how it is recovered when the setup is augmented allowing to trade these dynamically. This example is summarized in Figure \ref{fig:Ex1}. We consider a two period model with stock price process $S$ given by $S_0=S_1=0$ and $S_2\in\{-2,-1,1,2\}$. 
The American option process $\Phi$ is defined as $\Phi_1(\{S_1=0\})=1$, $\Phi_2(\{S_2\in\{-2,2\}\})=0$ and $\Phi_2(\{S_2\in\{-1,1\}\})=2$. The (pathwise) superhedging price of $\Phi$ can be easily computed and equals 2. A probability measure $\P$ on the space of four possible paths is uniquely described through a choice of $q_{i}=\P(S_2=i)\geq 0$ for $i\in\{-2,-1,1,2\}$ satisfying $q_2+q_1+q_{-1}+q_{-2}=1$. The martingale condition is equivalent to $2q_2+q_1-q_{-1}-2q_{-2}=0$. Note that there are only two stopping times greater than 0, $\tau_1=1$ and $\tau_2=2$, the market model price given as the double supremum over all stopping times $\tau$ and all martingale measures $\Q$ of $\E^\Q[\Phi_\tau]$ also equals 2 and the two prices agree. 
Suppose now that we add a {European option $g$} with a payoff $g=\id_{\{S_2=|1|\}}-1/2$ and initial price 0, which may be used as a  static hedging instrument. With $g$ and $S$, the superhedging price of $\Phi$ drops to $3/2$ (e.g.\ keep 3/2 in cash and buy one option $g$). In presence of $g$, we need to impose a calibration constraint on martingale measures: $q_1+q_{-1}=1/2$. Thus, any calibrated martingale measure can be expressed by $(q_2,q_1, q_{-1},q_{-2})=(q, 3/4-2q, 2q-1/4,1/2-q)$ with $q\in(1/8,3/8)$, and the market model price equals 1. We therefore see that adding a statically traded option breaks the pricing--hedging duality.
\tikzstyle{level 1}=[level distance=2.5cm, sibling distance=2cm]
\tikzstyle{level 2}=[level distance=3cm, sibling distance=1cm]
\tikzstyle{bag} = [text width=5.5em, text centered]
\tikzstyle{bagg} = [text width=3em, text centered]
\tikzstyle{end} = [circle, minimum width=3pt,fill, inner sep=0pt]
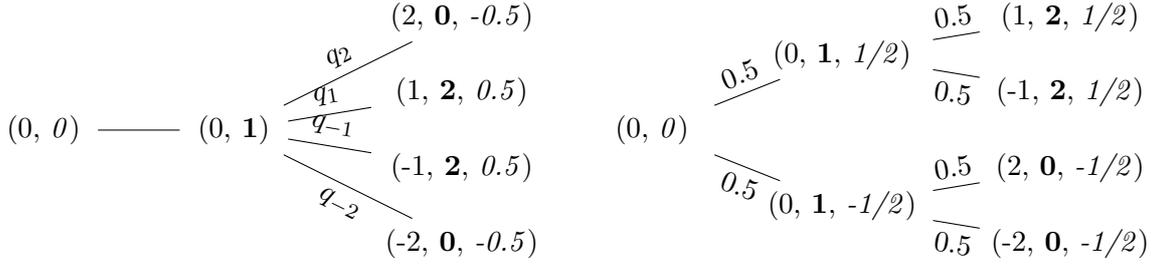
\begin{figure}
\begin{center}
\begin{tikzpicture}[grow=right, sloped]
\node[bagg] {({0}, \emph{0})}
    child {
        node[bagg] {({0}, \textbf{1}) }        
            child {
                node[bag]{({-2}, \textbf{0}, \emph{-0.5})}
                edge from parent
                node[below] {$q_{-2}$}
            }
            child {
                node[bag]{({-1}, \textbf{2}, \emph{0.5})}
                edge from parent
                node[above] {$q_{-1}$}
            }
            child {
                node[bag]{({1}, \textbf{2}, \emph{0.5})}
                edge from parent
                node[above] {$q_{1}$}
            }
            child {
                node[bag] {({2}, \textbf{0}, \emph{-0.5})}
                edge from parent
                node[above] {$q_{2}$}
            }
            edge from parent 
    };
\end{tikzpicture}
\begin{tikzpicture}[grow=right, sloped]
\node[bag] {({0}, \emph{0})}
    child {
        node[bag] {({0}, \textbf{1}, \emph{-1/2}) }        
            child {
                node[bag]{({-2}, \textbf{0},  \emph{-1/2})}
                edge from parent
                node[below] {0.5}
            }
            child {
                node[bag] {({2}, \textbf{0}, \emph{-1/2})}
                edge from parent
                node[above] {0.5}
            }
            edge from parent
            node[below] {0.5} 
    }
child {
        node[bag] {({0}, \textbf{1}, \emph{1/2})}   
            child {
                node[bag]{({-1}, \textbf{2}, \emph{1/2})}
                edge from parent
                node[below] {0.5}
            }
            child {
                node[bag]{({1}, \textbf{2}, \emph{1/2})}
                edge from parent
                node[above] {0.5}
            }
            edge from parent 
            node[above] {0.5}
    };
\end{tikzpicture}
\end{center}
\caption{Prices of the stock are written in regular font, prices of the American option in bold and prices of European option in italic. Model without dynamic trading in the option $g$ is on the left. Model with dynamic trading in $g$ which allows to recover the duality is on the right. \label{fig:Ex1}}
\end{figure}

Let us now show that the duality is recovered when we allow dynamic trading in $g$. We can model this through a process $Y=(Y_t:t=0,1,2)$ given  by $Y_2=g$, $Y_1=1/2$ on $\{S_2=|1|\}$, $Y_1=-1/2$ on $\{S_2=|2|\}$ and $Y_0=0$. Note that there exists a (unique) measure $\Q$ such that both $S$ and $Y$ are martingales w.r.t their joint natural filtration $\Fh$ and in particular $\Q$ is calibrated, $\E^\Q[g]=0$. The filtration $\Fh$ is richer than the natural filtration of $S$ alone and allows for an additional stopping time $\tau=\id_{\{Y_1=-1/2\}}+2\id_{\{Y_1=1/2\}}$ and the duality is recovered. 

\end{Example}


\section{Pricing--hedging duality for American options}
\label{duality}

We present in this section general results which explain when and why the pricing--hedging duality for American options holds. 
We work in a general discrete time setup which we now present. 
Let $(\Omega, \Fc)$ be a measurable space and $\ff:=(\Fc_k)_{k=0,1,...,N}$ be a filtration $\ff:=(\Fc_k)_{k=0,1,...,N}$, where $N\in \mathbb N$ is the time horizon.  
We denote by $\mathcal H$ the class of $\ff$-predictable processes. 
We denote by $\frak P(\Om)$ the set of all probability measures on $(\Om, \Fc)$. We consider a subset $\Pc\subset \frak P(\Om)$. 
We will say that a given property holds $\Pc$-quasi surely if it holds $\P$-almost surely for all $\P\in\Pc$. 
We will refer to a set from $\Fc$ to be a $\Pc$-polar set if it is a null set w.r.t. all $\P\in\Pc$. We will write $\Q\lll \Pc$ if there exists $\P\in\Pc$ such that $\Q\ll \P$.
Given a random variable $\xi$ and a sub-$\sigma$-field $\Gc \subset \Fc$, we define the conditional expectation $\E^{\P}[ \xi | \Gc] := \E^{\P}[ \xi^+ | \Gc] - \E^{\P} [ \xi^- | \Gc]$ with the convention $\infty -\infty = -\infty$, where $\xi^+ := \xi \vee 0$ and $\xi^- := - (\xi \wedge 0)$.
We consider a market with no transaction costs and with financial assets, some which are dynamically traded and some which are only statically traded. 
The former are modeled by an adapted $\R^d$-valued process $S$ with $d\in \mathbb N$. We think of the latter as European options which are traded at time $t=0$ but not necessarily at future times. 
We let $g=(g^\no)_{\no\in \No}$, where $\No$ is a set of an arbitrary cardinality, be the vector of their payoffs which are assumed $\Fc$-measurable and $\R$-valued. Up to a constant shift of the payoffs, we may, without loss of generality, assume that all options $g^\no$ have zero initial prices. Denote by $\Hc$ the set of all $\F$-predictable $\R^d$-valued processes, and by
$\hcal=\{h\in \R^\No: \exists \textrm{  finite subset $\beta\subset \No$ s.t. $h^\no=0$  $\forall \no\notin\beta$}\}$. A self--financing strategy trades dynamically in $S$ and statically in finitely many of $g^\no$, $\no \in \No$ and hence corresponds to a choice of $H\in \Hc$ and $h\in \hcal$. Its associated final payoff is given by
\be
\label{self}
(H \is S)_N + hg=\sum_{j=1}^d\sum_{k=1}^N H^j_k \Delta S^j_k + \sum_{\no\in\No} h^{\no} g^{\no},
\ee
where $\Delta S^j_k = S^j_{k}-S^j_{k-1}$. 
Having defined the trading strategy, we can consider the superhedging price of an option which pays off $\xi$ at time $N$:   
\begin{align}\label{eq:def_pi_E}
\pi^E_g(\xi):=\inf\{x: &\exists \ (H, h)\in \Hc\times\hcal \ \textrm{ s.t. } \ x+(H\is S)_N+hg\geq \xi \ \mathcal P\textrm{-q.s.} \}.
\end{align}
The inequality is required to hold $\Pc$-q.s., i.e.\ it holds $\P$-a.s.\ for any $\P\in \Pc$. In particular, if $\Pc = \Bf(\Om)$ is the set of all probability measures on $\Fc$ and $\{\omega\} \in \Fc$ for all $\omega \in \Omega$, then the superreplication in \eqref{eq:def_pi_E} is pathwise on $\Omega$. 

To formulate a duality relationship, we need the dual elements given by rational pricing rules, or martingale measures:
\be
\mathcal M&=&\{\mathbb Q\in \Bf(\Omega): \Q \lll \Pc \textrm{ and } \E^\Q[\Delta S_k| \Fc_{k-1}]=0, \, \forall k=1,...,N \}\nonumber\\
\label{eq:mcg}
\mathcal M_g&=&\{\mathbb Q\in \mathcal M: \E^\Q[g^\no]=0, ~ \forall \no\in \No \}.
\ee
\begin{Definition}
	Let $\fcl$ be a given class of functions defined on $\Omega$,
	we say that the \emph{(European) pricing--hedging duality holds for the class $\fcl$}
	if $\Mc_g\neq \emptyset$ and 
	\be\label{eq:PHduality}
		\pi^E_g(\xi)=\sup_{\Q\in \Mc_g}\E^\Q[\xi],\quad \xi\in\fcl,
	\ee
\end{Definition}

\begin{Remark}
\label{rem:weak}
Note that the inequality $"\geq"$ in \eqref{eq:PHduality}, called weak pricing--hedging duality, holds automatically from the definition of $\Mc_g$ in \eqref{eq:mcg}. 
\end{Remark}
A number of papers, including \cite{BN13,BFHMO16}, proved that the above pricing--hedging duality \eqref{eq:PHduality} holds under various further specifications and restrictions on $\Omega$, $\F$, $\Pc$ and $\fcl$, including in particular an appropriate no--arbitrage condition. We take the above duality for granted here and our aim is to study an analogous duality for American options. We work first in the general setup described above without specifying $\F$ or $\fcl$ since our results will apply to any such further specification. Further, many abstract results in this section also extend to other setups, e.g.\ to trading in continuous time.

\subsection{Superhedging of American options}
\label{superhedging}

An American option may be exercised at any time $k \in \T := \{1, \cdots, N\}$ (without loss of generality we exclude exercise at time $0$). It is described by its payoff function $\Phi = (\Phi_k)_{1 \le k \le N}$, where $\Phi_k: \Om \to \R$ belongs to $\fcl$ and is the payoff, delivered at time $N$, if the option is exercised at time $k$. Usually $\Phi_k$ is taken to be $\Fc_k$-measurable but here we only assume $\Phi_k$ to be $\Fc$-measurable for greater generality which includes, e.g.\ the case of a portfolio containing a mix of American and European options. 
We note that when hedging our exposure to an American option, we should be allowed to adjust our strategy in response to an early exercise. In consequence, the superhedging cost of the American option $\Phi$ using semi--static strategies is given by
\begin{align*}
\pi^A_g(\Phi)=\inf \big\{x: &\exists (H^1,...,H^N)\in \mathcal H^N \ \textrm{s.t.} \ H^j_i=H^k_i \ \forall 1\leq i\leq j\leq k\leq N\ \textrm{and} \ h\in \mathbb \hcal\\
&\textrm{satisfying} \ x+(H^k\is S)_N+hg\geq \Phi_k \ \forall k=1,...,N \ \mathcal P\textrm{-q.s.} \big\}
\end{align*}

Classically, pricing of an American option is recast as an optimal stopping problem and a natural extension of \eqref{eq:PHduality} would be
\be\label{eq:PHAduality}
\pi^A_g(\Phi)\stackrel{?}{=}\sup_{\Q\in\Mc_g}\sup_{\tau\in \Tc(\F)}\E^\Q[\Phi_\tau],
\ee
where $\Tc(\F)$ denotes the set of $\F$-stopping times. However, as illustrated with the simple example in the Introduction, this duality may fail. The ``numerical" reason is that the RHS in \eqref{eq:PHduality} may be too small since the set $\Mc_g$ is too small. Our aim here is to understand fundamental reasons why the duality fails and hence discuss how and why the right hand side should be modified to obtain equality in \eqref{eq:PHAduality}.

\subsection{American option is a European option on an enlarged space}
\label{ameu}

The first key idea of this paper offers a generic enlargement of the underlying probability space which turns all American options into European options. 
Depending on the particular setup, it may take more or less effort to establish \eqref{eq:PHduality} for the enlarged space but this shifts the difficulty back to the better understood and well studied case of European options. 
Our reformulation technique --- from an American to European option --- can be easily extended to other contexts, such as the continuous time case. 
The enlargement of space is based on construction of random times, previously used e.g. in \cite{JeanblancSong1, JeanblancSong2} to study the existence of random times with a given survival probability and in \cite{ET} to study a general optimal control/stopping problem, and in \cite{GTT2}, \cite{KTT} to study the optimal Skorokhod embedding problem. 

Let $\mathbb T:=\{1,..., N\}$ and introduce the probability space $\Omb:=\Omega \x \mathbb T$ with the canonical time $T: \Omb \to \T$ given by $T(\omb):=\theta$, where $\omb:=(\om, \theta)$, the filtration $\Fbb:=(\Fcb_k)_{k=0,1,...,N}$ with $\Fcb_k=\Fc_k\otimes \vartheta_k$ and $\vartheta_k=\sigma(T\land (k+1))$, and the $\sigma$-field $\Fcb=\Fc\otimes \vartheta_N$. 
By definition, $T$ is an $\Fbb$-stopping time. We denote bt $\Hcb$ the class of $\Fbb$-predictable processes and extend naturally the definition of $S$ and $g^\no$ from $\Om$ to $\Omb$ as $S(\omb)=S(\om)$ and $g^\no(\omb)=g^\no(\om)$ for $\omb=(\om, \theta)\in \Omb$.
We let $\fclb$ be the class of random variables $\xib:\Omb\to \R$ such that $\xib(\cdot,k)\in \fcl$ for all $k\in \T$ and we let $\pib^E_g(\bar{\xi})$ denote the superreplication cost of $\xib$. 
We may, and will, identify $\fclb$ with $\fcl^N$ via $\xib(\omb)=\Phi_\theta(\om)$.
Finally, we introduce
\be \label{eq:mcbg}
\Pcb&=&\{\Pb \in \frak P(\Omb): \Pb_{\vert\Om}\in \Pc\},\nonumber\\
\Mcb&=&\{\Qb \in \frak P(\Omb)  : \Qb\lll \Pcb \textrm{ and } \ \E^{\Qb}[\Delta S_k|\Fcb_{k-1}]=0 \ \forall k\in\T\}, \nonumber\\
\Mcb_g&=&\{\Qb \in\Mcb : \E^{\Qb}[g^\no]=0 \ \forall \no\in \No\}
\ee

\begin{Theorem}
\label{thm:American_European}
For any $\Phi\in \fcl^N = \fclb$ we have 
\be\label{eq:Amer_as_Euro}
	\pi^A_g(\Phi)=\pib^E_g(\Phi)
	:=
	\inf\{x: \exists \ (\Hb, h)\in \Hcb \times\hcal \ \textrm{ s.t. } \ x+(\Hb\is S)_N+hg\geq \xi \  \Pcb\textrm{-q.s.} \}.
\ee
In particular, if the European pricing--hedging duality on $\Omb$ holds for $\Phi$ then 
\be\label{eq:PHAduality_Euro}
\pi^A_g(\Phi)=\pib^E_g(\Phi)=\sup_{\Qb\in \Mcb_g}\E^{\Qb}[\Phi].
\ee
\end{Theorem}
%
\proof
First note that 
$$ \Hcb=\{\Hb=(\Hb(\cdot,1),\ldots,\Hb(\cdot,N))\in \Hc^N: \Hb_i(\cdot,j)=\Hb_i(\cdot,k) \ \forall 1\leq i\leq j\leq k\leq N
\}
$$
and hence the dynamic strategies used for superhedging in $\pi^A_g$ and in $\pib^E_g$ are the same. The equality now follows observing that  a set $\Gamma\in \Fcb_N$ is a $\Pcb$-polar if and only if its $k$-sections $\Gamma_k=\{\om: (\om, k)\in \Gamma\}$ are $\Pc$-polar  for all $k\in \T$. Indeed, for one implication assume that $\Pb(\Gamma)=0$ for each $\Pb\in \Pcb$. For arbitrary $\P\in \Pc$ and $k\in\T$ we can define $\Pb=\P\otimes \delta_k$ which belongs to $\Pcb$ and hence $\P(\Gamma_k)=0$.\\
To show the reverse implication, assume that $\P(\Gamma_k)=0$ for each $\P\in \Pc$ and $k\in\T$. Observe that, for any $\Pb\in\Pcb$
\begin{align*}
\Pb(\Gamma)=\sum_{k\in \T}\Pb(\Gamma_k\times \{k\})\leq\sum_{k\in \T}\Pb_{\vert\Om}(\Gamma_k)=0 
\end{align*}
since $\Pb_{\vert\Om}\in \Pc$. This completes the proof.
\qed
\begin{Remark}
\label{rem:g}
If the pricing--hedging duality holds w.r.t. filtration $\F$, 
then it holds as well for any filtration $\H\supset\F$ such that $\H$ 
and $\F$ only differ up to $\Mc_g$-polar sets. 
This change does not affect set $\Mc_g$ and may only decrease the superhedging cost as one has more trading strategies available. The duality is not affected by Remark \ref{rem:weak}. 
\end{Remark}
\begin{Remark}
\label{rem:weak_relaxed_formulation}
We note that the set $\Mcb_g$ in \eqref{eq:PHAduality_Euro}, or its projection on $\Omega$, is potentially much larger than $\Mc_g$. Indeed, instead of stopping times relative to $\F$, it allows us to consider any \emph{random} time which can be made into a stopping time under some calibrated martingale measure. We can rephrase this as saying that $\Mcb_g$ is equivalent to a \emph{weak formulation} of the initial problem at the r.h.s. of \eqref{eq:PHAduality}. To make this precise, let us call 
a weak stopping term $\alpha$ a collection
	$$
		\alpha ~=~ 
		\big( 
			\Om^{\alpha}, ~\Fc^{\alpha}, ~\Q^{\alpha}, \F^{\alpha} = (\Fc^{\alpha}_k)_{0 \le k \le N}, 
			~(S^{\alpha}_k)_{0 \le k \le N},
			~(g^{\no,\alpha})_{\no\in \No},~(\Phi_k^\alpha)_{k\in \T},
			~\tau^{\alpha}
		\big)
	$$
with $ \big( \Om^{\alpha}, ~\Fc^{\alpha}, ~\Q^{\alpha}, \F^{\alpha} \big)$ a filtered probability space, $\tau^{\alpha}$ a $\T$-valued $\F^\alpha$-stopping time, an $\R^d$-valued $(\Q^{\alpha}, \F^{\alpha})$-martingale $S^{\alpha}$ and a collection of random variables $g^{\no,\alpha},\Phi^\alpha_k$, and such that there is a measurable surjective mapping $\mathtt{i}_\alpha:\Om^\alpha\to \Om$ with $\Q=\Q^\alpha\circ \mathtt{i}_\alpha^{-1}\in \Mc$ and $\mathtt{i}_\alpha^{-1}(\Fc_k)\subset\Fc_k^\alpha$, $\mathtt{i}_\alpha^{-1}(\Fc)\subset\Fc^\alpha$ and finally $\Law_{\Q^\alpha}(S^\alpha, g^\alpha,\Phi^\alpha)=\Law_\Q(S, g,\Phi)$. 
	Denote by $\Ac_g$ the collection of all weak stopping terms $\alpha$ such that $\E^{\Q^{\alpha}}\big[ g^{\no,\alpha} \big] = 0$ for each $\no \in \No$. It follows that any $\alpha\in \Ac_g$ induces a probability measure $\Qb\in\Mcb_g$ and $\E^{\Q^{\alpha}}[\Phi^\alpha_{\tau^{\alpha}}]=\E^{\Qb}[\Phi]$. Reciprocally, any $\Qb\in\Mcb_g$, together with the space $(\Omb, \Fcb, \Fbb)$ and $(S, g, \Phi)$, provides a weak stopping term in $\Ac_g$. In consequence, 
		$$
			\sup_{\alpha \in \Ac_g} \!\E^{\Q^{\alpha}}\!\big[ \Phi_{\tau^{\alpha}} \big]
			=\!\!
			\sup_{\Qb \in \Mcb_g}\! \E^{\Qb}[ \Phi \big].
		$$
%
In summary, and similarly to number of other contexts, see the introduction in \cite{PhamZhang}, the weak formulation (and not the strong one) offers the right framework to compute the value of the problem. In fact, the set $\Mcb_g$ is large enough to make the problem static, or European, again. However, while it offers a solution and a corrected version of \eqref{eq:PHAduality}, it does not offer a fundamental insight into why \eqref{eq:PHAduality} may fail and what is the \emph{minimal} way of enlarging the objects on the RHS thereof to preserve the equality. These questions are addressed in the subsequent section.
\end{Remark}

\begin{Remark}
\label{com:hob}
\cite{Neuberger} and \cite{HobsonNeuberger} studied the same superhedging problem 
	in a Markovian setting, where the underlying process $S$ takes value in a discrete lattice $\Xc$.
	By considering the weak formulation (which is equivalent to our formulation, as shown in Remark \ref{rem:weak_relaxed_formulation} above), they obtain similar duality results.
	Moreover, they only consider $\Phi_k=\phi(S_k)$ where $\phi: \R^d\to\R$. Then the authors show that in the optimization problem $\sup_{\Qb\in \Mcb_g}\E^{\Qb}[\Phi]$ given in \eqref{eq:Amer_as_Euro} one may consider only Markovian martingale measures. The primal and the dual problem then turn to be linear programming problems under linear constraints, which can be solved numerically. Their arguments have also been extended to a more general context, where $S$ takes value in $\R_+$.
	Comparing to \cite{Neuberger, HobsonNeuberger}, our idea to weak formulation is very similar to theirs. However our setting is much more general and, when considering specific setups in Sections \ref{sec:main} and \ref{sec:mot}, we rely on entirely different arguments to prove the duality.
\end{Remark}

\subsection{The loss and recovery of the dynamic programming principle and the natural duality for American options}
\label{sec:dpp}

The classical pricing of American options, on which the duality in \eqref{eq:PHAduality} was modelled, relies on optimal stopping techniques which subsume a certain dynamic consistency, or a dynamic programming principle, as explained below. Our second key observation in this paper is that if the pricing--hedging duality \eqref{eq:PHAduality} for American options fails it is because the introduction of static trading of European options $g$ at time $t=0$ destroys the dynamic programming approach, the Bellman optimality principle. Indeed, $\pi^E_g(\xi)$ will typically be lower than the superhedging price at time $t=0$ of the capital needed at time $t=1$ to superhedge from thereon. To reinstate such dynamic consistency, we need to enlarge the model and consider dynamic trading in options in $g$. This will generate a richer filtration than $\F$ and one which will carry enough stopping times to obtain the correct natural duality in the spirit of \eqref{eq:PHAduality}. In particular, if $g=0$ (or equivalently $\Lambda = \emptyset$), then \eqref{eq:PHAduality} should hold. We now first prove this statement and then present the necessary extension when $g$ is non--trivial.

	Let $\fcl$ be a class of $\Fc$-measurable r.v.,
	we denote $\Ec(\xi) := \sup_{\Q\in \Mc}\E^{\Q}[\xi]$,
	and suppose that there is a family of operators $\Ec_k: \fcl\to \fcl$ for $k\in\{0,...,N-1\}$
	such that $\Ec_k(\xi)$ is $\Fc_k$-measurable for all $\xi \in \fcl$.
	We say that the family $(\Ec_k)$ provides a dynamic programming representation of $\Ec$ if 
\begin{equation}
\label{dpp_om}
\Ec(\xi)=\Ec_0\circ\Ec_1\circ...\circ\Ec_{N-1}\left(\xi\right),\ \ \forall \xi\in \fcl.
\end{equation}
The family $(\Ec_k)$ extends to the family $(\Ecb_k)$ for $k\in\{0,...,N-1\}$ defined for any $\Psi\in\fclb = \fcl^N$ by 
\begin{align}
\Ecb_0(\Psi)(\omb)&:=\Ec_0(\Psi(\cdot, 1))(\om), \quad \quad \quad\quad \quad \quad\quad \quad \quad ~\mbox{for all}~ \omb = (\om, \theta),
\label{eq:def_dpp0}\\
\Ecb_k(\Psi)(\omb)&:=
\begin{cases}
\Ec_k(\Psi(\cdot, \theta))(\omega)&\textrm{if }\theta<k\\
\Ec_k(\Psi(\cdot, k))(\omega)\lor \Ec_k(\Psi(\cdot, k+1))(\omega)&\textrm{if } \theta\geq k
\label{eq:def_dpp}
\end{cases}
, \quad \textrm{for $1\leq k<N$}.
\end{align}
Assume that $f\lor f'\in \fcl$ whenever $f,f'\in \fcl$, then $\Ecb$ maps functionals from $\fclb$ to $\fclb$.
Let us also introduce the $\Mc$-Snell envelope process of an American option $\Psi\in \fclb$ by
\be\label{eq:opbar}
\Ecb^k(\Psi):=\Ecb_k\circ ... \circ \Ecb_{N-1}(\Psi).
\ee
We say that the family $(\Ecb_k)$ provides a dynamic programming representation of $\Ecb(\Psi) := \sup_{\Q\in \Mcb}\E^{\Q}[\Psi]$ if 
\begin{equation}
\label{eq:dpp_omb}
\Ecb(\Psi) =\Ecb^0(\Psi),\ \ \forall \Psi\in \fclb.
\end{equation}
Typically we will consider $\Ec_k$ to be a supremum over conditional expectations w.r.t. $\Fc_k$, see Examples \ref{ex:dom} and \ref{ex:non} below. In these setups $\Ec_k$ automatically satisfies 
\begin{align}
\label{inequality2}
\sup_{\Qb\in \Mcb}\E^{\Qb}[\Psi]&\leq\Ecb^0(\Psi), \ \quad \Psi\in \fclb.
\end{align}
 
\begin{Theorem}
\label{dpp_gen_th}
Assume that $\No=\emptyset$, $\Ec_k$ satisfies  \eqref{dpp_om}, and that \eqref{inequality2} holds true, and that $f\lor f'\in \fcl$ for all $f,f'\in \fcl$.
Then, for all $\Phi\in\fcl^N = \fclb$,
\begin{equation}\label{eq:dpp_implied}
\sup_{\Qb\in \Mcb}\E^{\Qb}[\Phi]=\sup_{\Q\in \Mc}\sup_{\tau\in \Tc(\F)}\E^{\Q}[\Phi_\tau].
\end{equation}
If, further, the European pricing--hedging duality holds on $\Omb$ for the class $\fclb$, then
\begin{equation*}
\pi^A(\Phi)=\sup_{\Q\in \Mc}\sup_{\tau\in \Tc(\F)}\E^{\Q}[\Phi_\tau].
\end{equation*}
\end{Theorem}
The second assertion follows instantly from the first one and Theorem \ref{thm:American_European} while the first one follows from Proposition \ref{prop:dpp_Euro_American} below which asserts that \eqref{dpp_om} and \eqref{inequality2} imply an analogue consistency on $\Omb$ \eqref{eq:dpp_omb}. 
This also allows us to identify the optimal stopping time on the RHS of \eqref{eq:dpp_implied}. 
We have the following representation

\begin{Proposition}\label{prop:dpp_Euro_American}
Assume that $\No=\emptyset$ and $f\lor f'\in \fcl$ for all $f,f'\in \fcl$.
Then the dynamic programming representation \eqref{eq:dpp_omb} holds if and only if \eqref{dpp_om} and \eqref{inequality2} hold true.
Moreover, under condition \eqref{eq:dpp_omb}, the $\F$-stopping time
\begin{align}
\label{tau}
\tau^*(\omega)&:= \inf \left\{k\geq 1: \Ec_k\left(\Phi(\cdot, k)\right)(\omega)= \Ecb^{k}(\Phi)(\om, k)\right\}
\end{align}
provides the optimal exercise policy for $\Phi \in \fclb$:
\begin{equation}
\sup_{\Qb\in \Mcb}\E^{\Qb}[\Phi]=\sup_{\Q\in \Mc}\sup_{\tau\in \Tc(\F)}\E^{\Q}[\Phi_\tau]
=\sup_{\Q\in\Mc}\E^{\Q}\left[\Phi_{\tau^*}\right]=\Ecb^0(\Phi).
\end{equation}
\end{Proposition}

\begin{Remark} \label{rem:infinity_risky_assets}
	The proof of Proposition \ref{prop:dpp_Euro_American} will be provided in Section \ref{proofs1}. 
	The results in Theorem \ref{dpp_gen_th} and Proposition \ref{prop:dpp_Euro_American} are stated on $(\Om, \Fc)$, where there are only finitely many dynamic trading risky assets.
	However, its proof does not rely on the fact that the number of risky assets is finite,
	and the same results holds still true if there are infinitely number of dynamic trading risky assets.
\end{Remark}

	Next we give two examples of operators $(\Ec_k)_{k\leq N}$ satisfying \eqref{dpp_om} and \eqref{inequality2}, therefore, by Proposition \ref{prop:dpp_Euro_American}, also \eqref{eq:dpp_omb}.

\begin{Example}
\label{ex:dom}
Let $\Pc=\{\P^*\}$.  Then, taking $\fcl$ to be the set of all $\Fc$-measurable random variables and
\be 
\Ec_k(\xi)={\esup}_{\Q\in\Mc_g}\E^\Q[\xi|\Fc_k]
\ee
where $\esup$ is taken w.r.t $\P^*$, leads to a family of operators satisfying \eqref{dpp_om}, \eqref{inequality2}, and therefore also \eqref{eq:dpp_omb}. 
See the literature on dynamic coherent risk measures  for further discussion (e.g. \cite{AP} for an overview).
If in particular we assume that $\No=\emptyset$ then Theorem \ref{dpp_gen_th} recovers the classical superhedging theorem for American options (see e.g. \cite{Myneni}).
\end{Example}

\begin{Example}
\label{ex:non}
Let $(\Omega, d)$ be a Polish space, $\mathcal F$ its Borel $\sigma$-field and $\mathcal P$ a given set of probability measures on $(\Omega, \Fc)$. 
We are given a filtration $\ff:=(\Fc_k)_{k\leq N}$ such that 
$\Fc_0=\{\emptyset, \Omega\}$ and each $\sigma$-field $\Fc_k$ is countably generated.
For $k\in \T$ and $\om\in \Om$ denote by $\Mc^k(\om)$ the set of measures given by
$$\Mc^k(\om):=\{\Q\lll\Pc : \ \Q([\om]_{\Fc_k})=1\ \textrm{and}\ \E^\Q[\Delta S_n|\Fc_{n-1}]=0 \ \forall n\in\{k+1,..., N\} \}$$
where $[\om]_{\Fc_k}$ denotes the atom of $\Fc_k$ which contains $\om$, i.e.,
\be \label{eq:om_k}
	[\om]_{\Fc_k}=\bigcap_{F\in \Fc_k: \om\in F}F\ .
\ee
Note that $[\om]_{\Fc_k}\in \Fc_k$ since the latter is countably generated. In this setup we define  
\begin{equation*}
\Ec_k(\xi)(\om)=\sup_{\Q\in\Mc^{k}(\om)}\E^{\Q}\left[\xi\right].
\end{equation*}
If we furthermore assume that $\Ec_k(\xi)\in \fcl$ for any $\xi\in \fcl$ then the family $(\Ec_k)_{k\leq N}$ satisfies \eqref{inequality2} which will be proved in Proposition \ref{prop:non}.
We shall prove that under suitable assumptions on $(\Omega, \F,\Pc)$ and $\fcl$ also \eqref{dpp_om} holds for this family. 
This holds in particular in the setup of \cite{BN13} as shown therein, see (4.12) in \cite{BN13}.
\end{Example}

Let us consider the case with statically traded options: $\No\neq\emptyset$. We saw in Example \ref{ex:intro} that this can break down the dynamic consistency as the universe of traded assets differs at time $t=0$ and times $t\geq 1$. To remedy this, one has to embed the market into a fictitious larger one where both $S$ and all the options $g^\no$, $\no\in\No$, are traded dynamically. 
Let us consider a larger probability space $(\Omh,\Fch)$ which satisfies the following properties. 
First, there exists a surjective mapping $\mathtt{i}: \Omh\to \Om$, and it defines a natural extension of $S$, $g$ and $\Phi$ by $S(\omh)=S(\mathtt{i}(\omh))$, $g^\no(\omh)=g^\no(\mathtt{i}(\omh))$ for all $\no\in\No$ and $\Phi(\omh)=\Phi(\mathtt{i}(\omh))$.
Second there exists a family of processes $\Y=(\Y^\no)_{\no\in \No}$ on $\Omh$,  such that $\Y^\no_0=0$ and $\Y^\no_N(\omh)=g^\no(\mathtt{i}(\omh))$.
Let us denote by $\Sh=(S, \Y)$ which will now correspond to dynamically traded assets. 
We suppose that there is a filtration $\Fh:=(\Fch)_{k=0,1,...,N}$ such that $\mathtt{i}^{-1}(\Fc_k)\subset\Fch_k$ and $\Sh$ is $\Fh$-adapted, and let $\Hch$ be the set of $\Fh$-predictable processes. 
Finally, we consider the following sets of probability measures
\begin{align*}
\Pch&:=\{\Pbh\in \frak P(\Omh): \Pbh\circ\mathtt{i}^{-1}\in \Pc\ \}\\
\Mch&:=\{ \Qh\lll \Pch: \Sh=(S,\Y) \ \ \textrm{is an $(\Qh, \Fh)$-martingale}\}.
\end{align*}
Observe that the martingale measures in $\Mch$ are by definition calibrated to market prices of options in $g$.
We furthermore assume that the mapping $\mathtt{I}: \Mch\to\Mc_g$ defined by $\mathtt{I}(\Qh)=\Qh\circ\mathtt{i}^{-1}$ is surjective.
The collection $(\Omh, \Fch, \Fh, \mathtt{i}, \Y)$ satisfying above properties is called a dynamic extension of 
$(\Om,\Fc,\F, \Pc,S,g)$, or in short that $\Omh$ is a dynamic extension of $\Om$. 
Note that for any dynamic extension $\Omh$ it holds 
$$\Law_{\Qh}(S, Y)=\Law_{\mathtt{I}(\Qh)}(S, Y^{\mathtt{I}(\Qh)})$$
where $Y^{\no, \mathtt{I}(\Qh)}=(\E^\Q[g^\no|\Fc_k])_{k\leq N}$.
For any $\Qh\in\Mch$ let $\mathtt{I}(\Qh)=\Qh\circ\mathtt{i}^{-1}\in\Mc_g$. And conversely, from a given $\Q\in\Mc_g$ we may recover its ``parent" measure $\Qh\in\Mch$. 
We consider a class of functions $\fclh$ on $\Omh$ and assume that $\fcl\subset\fclh$ in the sense that for $f\in\fcl$, $f(\omh):=f(\mathtt{i}(\omh))$ belongs to $\fclh$.
Then the correspondence between $\Mch$ and $\Mc_g$ yields to 
$$
\sup_{\Q\in\Mch}\E^\Q\left[\xi\right]
= \sup_{\Q\in\Mc_g}\E^\Q\left[\xi\right]\quad \textrm{for any $\xi\in \fcl$}.
$$
As introduced at the beginning of Section \ref{ameu}, one can apply the enlargement techniques on space $\Omh$ to obtain $\Mct$, and one has a similar equality:
\be 
\label{eq:corr}
\sup_{\Q\in\Mct}\E^\Q\left[\Phi\right]
= \sup_{\Q\in\Mcb_g}\E^\Q\left[\Phi\right]\quad\textrm{for any $\Phi\in \fclb$}.
\ee
Since in the dynamic extension $(\Omh, \Fch, \Fh, \mathtt{i}, \Y)$ of $(\Om,\Fc,\F, \Pc,S,g)$ we allow to trade dynamically in $\Sh=(S, \Y)$ let us introduce the class of trading strategies $\Hh$ which is a set of $\Fh$-predictable $\R^{\Noh}$-valued processes which have only finitely many non--zero coordinates where $\Noh=\{(i,s): i\in\{1,..,d\}\}\cup\{(\no,y):\no\in\No\}$, i.e., 
\begin{align*}
\Hch=\{\Hh=&{(\Hh^{\noh}_k: {\noh\in\Noh})}_{k\leq N}: \  \textrm{$\Fh$-predictable $\R^{\Noh}$-valued process s.t.}\\
&\textrm{ $\exists$ finite subset $\Noh_0 \subset \Noh$ s.t. $\Hh^{\noh}_k=0, ~\forall k,~\forall \noh\notin\Noh_0$  }  \}.
\end{align*}
In consequence, a self--financing strategy corresponds to a choice of $\Hh\in \Hch$ and yields a final payoff of
\be
\label{self}
(\Hh \is \Sh)_N=\sum_{j=1}^d\sum_{k=1}^N \Hh^{(j,s)}_k \Delta S^j_k +\sum_{\no\in\No}\sum_{k=1}^N \Hh^{(\no,y)}_k \Delta Y^\no_k.
\ee
Note that appropriate choice of trading strategies ensures that the sums are finite.
The supehedging costs of a European option $\xih$ and an American option $\Psih=(\Psih_k)_{k\leq N}$ on $\Omh$ are given by 
\begin{align*}
\pih^E(\xih)=\inf\{x: &\exists \Hh\in \Hch \ \textrm{s.t.} \ \textrm{satisfying} \ x+(\Hh\is \Sh)_N\geq \xih\ \Pch\textrm{-q.s.} \},\\
\pih^A(\Psih)=\inf\{x: &\exists (\Hh^1,...,\Hh^N)\in \Hch^N \ \textrm{s.t.} \ \Hh^j_i=\Hh^k_i \ \forall \;1\leq i\leq j\leq k\leq N\ \\
&\textrm{satisfying} \ x+(\Hh^k\is \Sh)_N\geq \Psih_k \ \forall k=1,...,N \ \Pch\textrm{-q.s.} \}.
\end{align*}
\begin{Remark}
Clearly $\Fh$ is much richer than $\F$.
Beside the price process of the underlying, it captures all possible price processes of vanilla options. Therefore the inequality $\pih^A(\Phi)\leq \pi^A_g(\Phi)$ holds which follows by noting that buy--and--hold strategies are a special case of a dynamic trading strategy and $\Pc=\Pch\circ \mathtt{i}^{-1}$.
\end{Remark}

We can now apply Theorem \ref{dpp_gen_th} to the present setting:
\begin{Corollary}
\label{cor:dpp_options}
Let $(\Omh, \Fch, \Fh, \mathtt{i}, \Y)$ be the dynamic extension of $(\Om,\Fc,\F, \Pc,S,g)$ with operators $\Ech_k:\fclh\to\fclh$ satisfying \eqref{dpp_om} and \eqref{inequality2}.
Assume that the European pricing--hedging duality holds for the class $\fcl^N$ on $\Omb$. 
Then for all $\Phi\in\fcl^N$
\begin{equation}\label{eq:PHA_dpp}
\pi^A_g(\Phi)= \pih^A(\Phi)=\sup_{\Qh\in \Mch}\sup_{\tau\in \Tc(\Fh)}\E^{\Qh}\left [\Phi_\tau\right].
\end{equation}
\end{Corollary}
\proof
Note that $\pi^A_g\geq \pih^A$ since a buy--and--hold strategy is a special case of a dynamic trading strategy and $\Pc=\Pch\circ \mathtt{i}^{-1}$.
Using \eqref{eq:Amer_as_Euro} twice we obtain
$$\pib^E_g(\Phi)=\pi^A_g(\Phi)\geq \pih^A(\Phi)=
\pit^E\left(\Phi\right)\geq
\sup_{\Q\in\Mct}\E^\Q\left[\Phi\right]
= \sup_{\Q\in\Mcb_g}\E^\Q\left[\Phi\right],
$$
where the penultimate inequality always holds by Remark \ref{rem:weak} and last equality follows by \eqref{eq:corr}. The assumed pricing--hedging duality on $\Omb$ implies that we have equality throughout and we conclude by 
applying Theorem \ref{dpp_gen_th} (with Remark \ref{rem:infinity_risky_assets}) on $\Omh$ for the representation of $\sup_{\Q\in\Mct}\E^\Q\left[\Phi\right]$.
\qed

\begin{Remark}
Note that if pricing--hedging duality on $\Omb$ holds then trading vanilla options dynamically or statically makes no difference on superhedging cost.
\end{Remark}

\begin{Example}
\label{ex:product}
We give an example of a dynamic extension of $\Om$ in the case of finitely many statically traded options, i.e. we assume that $\No=\{1,...,e\}$ for some $e\in\mathbb N$.
Consider the probability space $\Omh=\Om\times \R^{(N-1)\times e}$. 
An element $\omh$ in $\Omh$ can be written as $\omh=(\om, \y)$ where $\y=(\y^1, ..., \y^e)\in \R^{(N-1)\times e}$ with $\y^i=(\y^i_1, ..., \y^i_{N-1})$. Define a mapping $\mathtt{i}:\Omh\to\Om$ by $\mathtt{i}(\omh)=\om$ which is clearly surjective.
We also introduce the process $\Y$ as $\Y_k(\omh)=\y_k=(\y^1_k, ..., \y^e_k)$ for $k\in\{1,...,N-1\}$, $\Y_0(\omh)=0$ and $\Y_N(\omh)=g(\omh)=g(\om)$. 
Let the filtration $\Fh:=(\Fch)_{k=0,1,...,N}$ be given by $\Fch_k=\Fc_k\otimes \Yc_k$, $\Yc_k=\sigma(\Y_n: n\leq k)$. 
In this case it also holds that $\mathtt{I}:\Mch\to \Mc_g$ given by $\mathtt{I}(\Qh)=\Qh\circ\mathtt{i}^{-1}=\Qh_{\vert\Om}$ is surjective.
In Section \ref{sec:main}, where the basic setup is taken from \cite{BN13} and hence \eqref{dpp_om} holds on $\Om$ as recalled above in Example \ref{ex:non}, we show that \eqref{dpp_om} also holds on $\Omh$.
\end{Example}

\begin{Remark}
Let us consider the two period ($N=2$) example of \cite{HobsonNeuberger2}, see Figure \ref{fig:hobson} below. 
For simplicity, we introduce only one statically traded option $g$ with payoff $\id_{\{S_2=8\}}$ at time $t=2$ and price $2/5$ at time $t=0$. This already destroys the pricing--hedging duality for the American option $\Phi$.
In \cite{HobsonNeuberger2}, the duality is recovered by considering a (calibrated) mixture of martingale measures. 
It is insightful to observe that their mixture model is nothing else but a martingale measure for an augmented setup with dynamic trading in $g$ which, following Corollary \ref{cor:dpp_options}, restores the dynamic programming and the pricing--hedging duality for American options. 
To show this, let $Y$ denote the price process of the option $g$: $Y_0=2/5$, $Y_2=g$. Figure \ref{fig:hobson} illustrates a martingale measure $\Q$ along with the intermediate prices $Y_1$ such that the processes $S$ and $Y$ are martingales. With $\tau=\id_{\{S_1=1, Y_1=0\}}+2\id_{\{S_1=1, Y_1=1/4\}\cup\{S_1=3\}}$ we find  $\E^\Q[\Phi_\tau]=18/5$ which is the super--hedging price and the duality is recovered. 
\tikzstyle{level 1}=[level distance=3.5cm, sibling distance=2.7cm]
\tikzstyle{level 2}=[level distance=3cm, sibling distance=1cm]

\tikzstyle{bag} = [text width=5em, text centered]
\tikzstyle{end} = [circle, minimum width=3pt,fill, inner sep=0pt]
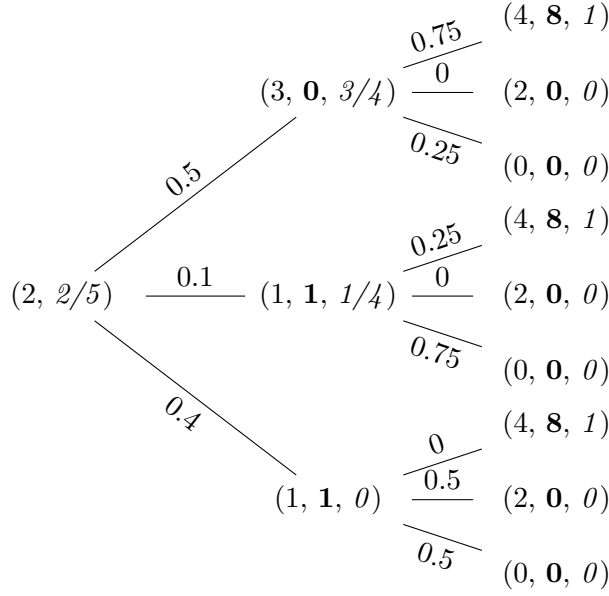
\begin{figure}
\label{fig:hobson}
\begin{center}
\begin{tikzpicture}[grow=right, sloped]
\node[bag] {({2}, \emph{2/5})}
    child {
        node[bag] {({1}, \textbf{1}, \emph{0})}        
            child {
                node[bag]{({0}, \textbf{0}, \emph{0})}
                edge from parent
                node[below] {0.5}
            }
            child {
                node[bag]{({2}, \textbf{0}, \emph{0})}
                edge from parent
                node[above] {0.5}
            }
            child {
                node[bag]{({4}, \textbf{8}, \emph{1})}
                edge from parent
                node[above] {0}
            }
            edge from parent 
            node[below] {0.4}
    }
child {
        node[bag] {({1}, \textbf{1}, \emph{1/4}) }        
            child {
                node[bag]{({0}, \textbf{0}, \emph{0})}
                edge from parent
                node[below] {0.75}
            }
            child {
                node[bag]{({2}, \textbf{0}, \emph{0})}
                edge from parent
                node[above] {0}
            }
            child {
                node[bag]{({4}, \textbf{8}, \emph{1})}
                edge from parent
                node[above] {0.25}
            }
            edge from parent 
            node[above] {0.1}
    }
    child {
        node[bag] {({3}, \textbf{0}, \emph{3/4})}        
            child {
                node[bag]{({0}, \textbf{0}, \emph{0})}
                edge from parent
                node[below] {0.25}
            }
            child {
                node[bag]{({2}, \textbf{0}, \emph{0})}
                edge from parent
                node[above] {0}
            }
            child {
                node[bag]{({4}, \textbf{8}, \emph{1})}
                edge from parent
                node[above] {0.75}
            }
            edge from parent 
            node[above] {0.5}
    };
\end{tikzpicture}
\end{center}
\caption{The model on $\Omh$ which corresponds to mixture model in \cite{HobsonNeuberger2} attaining the super--hedging price. Prices of the stock are written in regular font, prices of the American option in bold and prices of European option in italic.}
\end{figure}
\end{Remark}

\subsection{Pseudo--stopping times}
It follows from Theorem \ref{thm:American_European} that in general we expect to see
$$
\pi^A_g(\Phi)=\sup_{\Qb\in \Mcb_g}\E^{\Qb}[\Phi]\geq \sup_{\Q\in\Mc_g}\sup_{\tau\in \Tc(\F)}\E^\Q[\Phi_\tau],
$$
where the last inequality may be strict. We showed above that this is linked with the necessity to use random times beyond $\tau\in \Tc(\F)$. To conclude our general results, we explore this property from another angle and identify the subset(s) of $\Mcb_g$ which leads to equality in place of inequality above. Introduce 
	\be
		~~\Mce_g ~:=~
		\big\{\Qb \in \Bf(\Omb) \!\!\!\!\!&:&\!\!\!  \Qb \lll \Pcb,  ~\E^{\Qb}[g^\no] \!=\! 0, \;\;\no \in\No \;\;
~ S~\mbox{is an}~(\F, \Qb) \mbox{-martingale,}\nonumber\\
\E^{\Qb}[M_T] &=& \E^{\Qb}[M_0]
		~~\mbox{for all bounded $(\F, \Qb)$-martingales}~M
		\big\},
\ee
the set of measures which make $S$ an $\F$-martingale and $T$ an $\F$-pseudo--stopping time. These are natural since the martingale part of the Snell envelope can be stopped at the pseudo--stopping time with null expectation.
\begin{Proposition} 
\label{pseudo}
Assume that $\Mc_g\neq \emptyset$. Then
\be \label{eq:equiv_immersion}
\sup_{\Qb \in \Mce_g} \E^{\Qb}[\Phi]
=\sup_{\Q \in \Mc_g} \sup_{\tau \in \Tc(\F)} \E^{\Q} \big[ \Phi_{\tau}].
\ee
\end{Proposition}
\proof
Let $\Qb \in \Mce_g$ such that $\E^{\Qb}[|g^\no|] < \infty$ and $\E^{\Qb}[|\Phi_k|] < \infty$ for all $\no\in\No$ and $k=1, \cdots, N$.
We next consider the optimal stopping problem $\sup_{\tau \in \Tc(\F)} \E^{\Qb} \big[ \Phi_{\tau} \big]$.	
Define its Snell envelope $(Z_k)_{0 \le k \le N}$ by
	$$
		Z_k := \mathrm{esssup}_{\tau \in \Tc(\F), \tau \ge k} \E^{\Qb} \big[\Phi_{\tau} \big| \Fc_k \big],
	$$
which is an $(\F,\Qb)$-supermartingale. Its Doob--Meyer decomposition is given by
	$$
		Z_k = Z_0 + M_k - A_k, ~~~\mbox{where}~A=(A_k)_{0 \le k \le N}~\mbox{is an $\F$-predictable increasing process},
	$$	
	and $A_0 = M_0 = 0$.
	It follows that
	\be
		\E^{\Qb} \big[ \Phi \big]=\E^{\Qb} \big[ \Phi_T \big]
		~\le~
		\E^{\Qb} [Z_T]
		~\le~
		Z_0 + \E^{\Qb}[M_T]
		~=~
		Z_0.
	\ee
	We hence obtain that $\sup_{\Qb \in \Mce_g} \E^{\Qb}[\Phi] \le \sup_{\Q \in \Mc_g} \sup_{\tau \in \Tc(\F)} \E^{\Q} \big[ \Phi_{\tau}]$.
	Then \eqref{eq:equiv_immersion} holds since every stopping time $\tau \in \Tc(\F)$ is a pseudo-stopping time and hence the inverse inequality is trivial.
\qed
	
\begin{Remark}
The above allows us to see that it is not enough to use randomised stopping times to recover equality in \eqref{eq:PHAduality}. 
Such a time corresponds to an $\F$-adapted increasing process $V$ with $V_0=0$ and $V_N=1$. It may be seen as a distribution over all possible stopping times, in our setup a distribution $\eta$ on $\T$ s.t. $\eta(\{k\}):=\Delta V_k=V_k-V_{k-1}$ for each $k\in \T$. 
For any pseudo--stopping time $\tau$, the dual optional projection of the process $\id_{\Lbrack \tau, N\Rbrack}$ is a randomised stopping time. 
Conversely, for a given $V$, if we take a uniformly distributed random variable $\Theta$ independent from $V$, possibly enlarging probability space, then $\tau:=\inf\{t: V_t\geq \Theta\}$ is $\F$-pseudo--stopping time which generates $V$.
Let $\Rc$ be the set of such randomised stopping times. Then, from Proposition \ref{pseudo} and definition of dual optional projection, 
\begin{equation*}
\sup_{\Q \in \Mc_g} \sup_{\tau \in \Tc(\F)} \E^{\Q} \big[ \Phi_{\tau}]=\sup_{\Q \in \Mc_g} \sup_{V \in \Rc} \E^{\Q} \left[ \sum_k\Phi_{k}\Delta V_k\right].
\end{equation*}
\end{Remark}
\begin{Remark}
\cite{NY05} showed that under a progressive enlargement with pseudo--stopping time $\tau$ all martingales from the smaller filtration stopped at $\tau$ remain martingales in the larger filtration. 
One can relate this to a more restrictive situation, when all martingales from the smaller filtration remain martingales in the bigger filtration, which is called the immersion property in enlargement of filtration. Clearly each random time satisfying immersion property is a pseudo--stopping time. Thus, keeping the equality \eqref{eq:equiv_immersion} true, the pseudo--stopping time property in the definition of $\Mce_g$ above can be replace by a stronger condition characterizing the immersion property
		\be \label{eq:cond_immersion}
			\Qb[T > k | \Fc_n] = \Qb[ T >k | \Fc_k], 
			~~\mbox{for all}~
			0 \le k \le n \le N.
		\ee
See Section 3.1.2 of \cite{BJR} for the discrete time context of progressive enlargement of filtration and \cite{AL16} for connections between pseudo--stopping times, the immersion property and projections.
	\end{Remark}

\section{A detailed study of the non--dominated setup of\\ \cite{BN13}}
\label{sec:main}

In this section we work in the non--dominated setup introduced in \cite{BN13} which is a special case of Example \ref{ex:non}. 
We let $\Om_0 =\{\om_0\}$ be a singleton and $\Om_1$ be a Polish space. 
	For each $k \in \{1, \cdots, N\}$, we define $\Om_k := \{\om_0\}\times\Om_1^k$ as the $k$-fold Cartesian product.
	For each $k$, we denote by $\Gc_k := \Bc(\Om_k)$ and by
	$\Fc_k$ its universal completion.	
	In particular, we notice that $\Gc_0$ is trivial and we denote 
	$$
		\Om := \Om_N, ~~ \Gc:=\Gc_N~~\mbox{and}~~ \Fc := \Fc_N.
	$$
	We shall often see $\Fc_k$ and $\Gc_k$ as sub-$\sigma$-fields of $\Fc_N$, 
	and hence obtain two filtrations $\F = (\Fc_k)_{0 \le k \le N}$ and $\G = (\Gc_k)_{0 \le k \le N}$ on $\Om$.	
Recall that a subset of a Polish space $\Om$ is analytic if it is the image of a Borel subset of another Polish space under a Borel measurable mapping. 
	We take $\fcl$ to be the class of upper semianalytic functions $f: \Om \to \Rb := [-\infty, \infty]$, i.e.\ such that $\{\om \in \Om ~: f(\om) > c\}$ is analytic for all $c \in \R$.

The price process $S$ is a $\G$-adapted $\R^d$-valued process and the collection of options $g=(g^1, ..., g^e)$ is a $\G$-measurable $\R^e$-valued vector for $e\in \mathbb N$ (thus $\No=\{1,...,e\}$).

Let $k \in \{0, \cdots, N-1\}$ and $\om \in \Om_k$, 
	we are given a  non--empty convex set $\Pc_k(\om) \subseteq \Bf(\Om_1)$ of probability measures, which represents the set of all possible models for the $(k+1)$-th period, given state $\om$ at times $0, 1, \cdots, k$.
	We assume that for each $k$, 
	\be \label{eq:AnalyticGraph}
		\mbox{graph}(\Pc_k)
		:= 
		\{ (\om, \P): \om \in \Om_k, 
			\P \in \Pc_k(\om) 
		\} 
		~\subseteq~ \Om_k \times \Pc(\Om_1)
		~~ \text{is analytic.}
	\ee
	Given a universally measurable kernel $\P_k: \Om_k \to \Bf(\Om_1)$ for each $k \in \{ 0,1, \cdots, N-1 \}$, 
	we define a probability measure $\P$ on $\Om$ by Fubini's theorem:
	$$
 		\P(A) 
		:=
		\int_{\Om_1} \cdots \int_{\Om_1} \mathbf{1}_A (\om_1, \om_2 \cdots, \om_N) 
		\P_{N-1}(\om_1, \cdots, \om_{N-1}; d \om_N)  \cdots \P_0(d \om_1).
	$$
	We can then introduce the set $\Pc \subseteq \Bf(\Om)$ of possible models for the multi--period market up to time $N$:
	\be \label{eq:def_Pc}
		\Pc := 
		\big\{
			\P_0 \otimes \P_1 \otimes \cdots \otimes \P_{N-1} ~: \P_k(\cdot) \in \Pc_k(\cdot), k=0, 1, \cdots, N-1 
		\big\}.
	\ee
	Notice that the condition \eqref{eq:AnalyticGraph} ensures that $\Pc_k$ has always a universally measurable selector: $\P_k : \Om_k \rightarrow \Bf(\Om_1)$ such that $\P_k(\om) \in \Pc_k(\om)$ for all $\om \in \Om_k$.
	Then the set $\Pc$ defined in \eqref{eq:def_Pc} is nonempty.
We also denote by $\Mc^{k,k+1}(\om)$ the following set
$$\Mc^{k,k+1}(\om)=\{\Q\in \mathfrak P(\Om_1): \Q \lll \Pc_k(\om) \ \textrm{and} \ \ \E^{\delta_{\om} \otimes_k \Q}[ \Delta S_{k+1}] = 0\},$$
	where $\delta_{\om} \otimes_k \Q := \delta_{(\om_1, \cdots, \om_k)} \otimes \Q$ is a Borel probability measure on $\Om_{k+1} := \Om_k \x \Om_1$.

	The following notion of no--arbitrage $\NA(\Pc)$ has been introduced in \cite{BN13}.
$\NA(\Pc)$ holds if for all $(H, h) \in \Hc \times \mathbb{R}^e$
		$$
			(H \is S)_N + hg \geq 0 
			\quad \Pc\mbox{-q.s.}~
			\Longrightarrow~ 
			(H \is S)_N + hg = 0
			\quad \Pc\mbox{-q.s.}
		$$
Analogously, we will say that $\NA(\Pcb)$ holds if for all $(\Hb, h) \in \Hcb \times \mathbb{R}^e$
\be\label{nab}
			(\Hb \is S)_N + hg \geq 0\quad	\Pcb\mbox{-q.s.}~
			\Longrightarrow~ 
			(\Hb \is S)_N + hg = 0
			\quad \Pcb\mbox{-q.s.}
\ee
	Recall also $\Mc_g$ and $\Mcb_g$ have been defined in \eqref{eq:mcg} and \eqref{eq:mcbg}.
	As established in \cite{BN13}, the condition $\NA(\Pc)$ is equivalent to the statement that $\Pc$ and $\Mc_g$ have the same polar sets.
The following lemma extends that result to $\Omb$.
\begin{Lemma}
\label{lem:na}
$\NA(\Pc)$ $\Longleftrightarrow$ $\NA(\Pcb)$ $\Longleftrightarrow$ $\Pcb$ and $\Mcb_g$ have the same polar sets.
\end{Lemma}

\proof
The two conditions $\NA(\Pc)$ and $\NA(\Pcb)$ are equivalent by the same arguments as in proving \eqref{eq:Amer_as_Euro}.
It is enough to show that $\Pcb$ and $\Mcb_g$ have the same polar sets if and only if $\Pc$ and $\Mc_g$ have the same polar sets.
That boils down to proving that a set $\Gamma\in \Omb$ is an $\Mcb_g$ polar set if and only if  the $k$-section $\Gamma_k=\{\om: (\om, k)\in \Gamma\}$ is an $\Mc_g$ polar set for each $k\in \T$ which is analogous statement involving $\Pc$ and $\Pcb$ proved in the proof of Theorem \ref{thm:American_European}. 
\qed

\subsection{Duality on the enlarged space $\Omb$}

Our first main result is the following duality  under the no--arbitrage condition \eqref{nab}.

\begin{Theorem} \label{theo:main}
Let $\NA(\Pcb)$ hold true.
Then the set $\Mcb_g$ is nonempty, and, for any upper semianalytic $\Phi: \Omb \to \R$, one has
\be \label{eq:main}
\pib^E_g(\Phi) &=& \sup_{\Qb \in \Mcb_g} \E^{\Qb}[ \Phi \big],
\ee
and in particular the pricing-hedging duality \eqref{eq:PHAduality_Euro} holds. 
Moreover, there exists $(\Hb, h) \in \Hcb \x \R^e$ such that
$$\pib^E_g(\Phi) + (\Hb \is S)_N + h g \ge \Phi, ~~~\Pcb \mbox{-q.s.}
$$
\end{Theorem}
The proof is delegated to Section \ref{proofs2} and uses the following lemma.
Let us work with operators $\Ec_k$ introduced in Example \ref{ex:non}. Observe that 
$$\Ec_k\circ ... \circ \Ec_{N-1}(\xi)(\om)=\Ec_{k,k+1}\circ ... \circ \Ec_{N-1,N}(\xi)(\om),\quad \xi\in \fcl$$
where $\Ec_{k,k+1}(\xi)(\om)=\sup_{\Mc^{k,k+1}(\om)}\E^{\Q}\left[\xi\right]$. By Proposition \ref{prop:dpp_Euro_American}, (4.12) in \cite{BN13} and using that the maximum of upper semianalytic functions is still upper semianalytic we conclude that

\begin{Lemma} \label{lemma:super_hedging}
Consider the case $e=0$, i.e., $\No=\emptyset$.
Let $\Psi\in\fclb$. Then $\Ecb_k(\Psi)$ in \eqref{eq:def_dpp} is also upper semianalytic and
		$$
			\sup_{\Qb \in \Mcb} \E^{\Qb}[ \Psi]	
			~=~
			\Ecb^0(\Phi)
			:=
			\Ecb_0 \circ \cdots \circ \Ecb_{N-1} (\Psi).
		$$
\end{Lemma}

\subsection{Dynamic programming principle on $\Omh$}

Recall that the family of operators $(\Ec_k)$ on functionals on $\Om$ is defined in Example \ref{ex:non}, 
based on which one obtains a family of operators $(\Ech_k)$ on functionals on $\Omh$ as in Example \ref{ex:product}.

\begin{Theorem}
	Let $\xih: \Omh \to \R^N$ be an upper semianalytic functional. 
	Then $\Ect_k(\xih)$ is also upper semianalytic and
		$$
			\sup_{\Qh \in \Mct_g} \E^{\Qh}[ \xih]	
			~=~
			\Ect^0(\xih)
			:=
			\Ect_0 \circ \cdots \circ \Ect_{N-1} (\xih).
		$$
In particular, $\NA(\Pc)$ implies that \eqref{eq:PHA_dpp} holds.
\end{Theorem}
\proof
By Lemma \ref{lem:na} $\NA(\Pc)\iff\NA(\Pcb)$, and then by Theorem \ref{theo:main} the pricing--hedging duality on $\Omb$ in \eqref{eq:main} holds.
Then, by Corollary \ref{cor:dpp_options}, \eqref{eq:PHA_dpp} is implied by the dynamic programming principle on $\Omh$ for which it is enough to argue that the $\Mch$ satisfies the same analyticity property as $\Pc$, or $\Mc$,
	and the assertion would follow by (4.12) in \cite{BN13}.
	
	Given $k =0, \cdots, N-1$ and $\om \in \Om$,
similar to \eqref{eq:def_Pc}, 
we define
$$
		\Mc'^{,k}(\om) ~:=~ \{\Q := \om \otimes_k \Q_{k+1} \otimes \cdots \otimes \Q_{N-1} ~: \Q_i \in \Mc^{i,i+1}(\om)
		 \}.
	$$
The collection $\Mc'^{,k}(\om)$ induces a collection $\Mch'^{,k}(\om,y)$ on $\Omh$ by
	$$
		\Mch'^{,k}(\om, y) ~:=~ \{\Qh := \Q \circ (X, Y^{\Q})^{-1}, ~\Q \in \Mc'^{,k}(\om),~ Y^{\Q}_k = y\}.
	$$
	By considering its marginal law at time $k+1$, we define
	$$
		\Mch^k(\om,y) ~:=~\{ \Q \circ (X_{k+1}, Y^{\Q}_{k+1})^{-1}, ~\Q \in \Mc'^{,k}(\om),~ Y^{\Q}_k = y \}.
	$$
	We claim that
	\be \label{eq:claim_analytic}
		\mbox{the graph}~
		\{ (\om, y, \Qh) ~: \Qh \in \Mch^k(\om, y)\}
		~\mbox{is analytic.}
	\ee
	Then the problem reduces to the same context as in Lemma \ref{lemma:super_hedging} and one obtains immediately the dynamic programming representation for $\sup_{\Qh \in \Mct_g} \E^{\Qh}[\cdot]$ as stated.
	
	To conclude the proof, it is enough to prove the claim \eqref{eq:claim_analytic}.
	First, since the graph $\{(\om, \Q) ~: \Q \in \Mc^{k,k+1}(\om)\}$ is analytic,
	it follows by Theorem 2 in \cite{DellacherieMaisonJeu} 
	that the graph
	$\{(\om, \Q) ~:= \Mc'^{,k}(\om) \}$ is also analytic. 
	(Notice that the results in \cite{DellacherieMaisonJeu} is given in a Markovian context, 
	by considering the whole path, we can easily reduce our problem to his Markovian context.)
	Next, by Lemma 3.1 of \cite{NN14}, we can choose a version of a family $(Y^{\Q})_{\Q \in \Bc}$ such that
	$(\om, \Q) \mapsto Y^{\Q}(\om)$ is Borel measurable.
	It follows that the graph
	$\{(\om, y, \Q) ~: \Q \in \Mch'^{,k}(\om,y)\}$ is analytic,
	thus 
	$\{(\om, y, \Q) ~: \Q \in \Mch'^{,k}(\om,y)\}$ is also analytic.
	\qed
	
\subsection{Comparison with \cite{Bay} and \cite{BayZhou}}
\label{com:bay}
	
	In \cite{Bay}, the authors considered the same superhedging problem $\pi^A_g(\Phi)$ with finite set $\No=\{1,...,e\}$, 
	and established the duality
	\be \label{eq:dualBay}
		\pi^A_g(\Phi) = \inf_{h \in \R^e} \sup_{\tau \in \Tc(\F)} \sup_{\Q \in \Mc_0} \E^{\Q}[ \Phi_{\tau} - hg],
	\ee
	under some regularity conditions (see Proposition 3.1 in \cite{Bay}).
	Our duality in Theorem \ref{theo:main} is more general and more complete,
	and moreover, together with Lemma \ref{lemma:super_hedging}, it induces the above duality \eqref{eq:dualBay}.
	In exchange, \cite{Bay} also studied another subhedging problem $\sup_{\tau \in \Tc(\F)} \inf_{\Q \in \Mc} \E^{\Q}[\Phi_{\tau}]$ which we do not consider here.
	
	More recently, \cite{BayZhou}  consider the ``randomized'' stopping times, 
	and obtain a more complete duality for $\pi^A_g(\Phi)$.
	The dual formulations in \cite{BayZhou} and in our results are more or less in the same spirit (as in \cite{Neuberger, HobsonNeuberger}).
	Nevertheless, the duality in \cite{BayZhou} is established under strong integrability conditions and an abstract condition which is checked under regularity conditions (see their Assumption 2.1 and Remark 2.1).
	In particular, when $\Pc$ is the class of all probability measures on $\Om$, the integrability condition in their Assumption 2.1 is equivalent to say that $\Phi_k$ and $g^i$ are all uniformly bounded.
	In our paper, we only assume that $g^i$ are Borel measurable, $\Phi_k$ are upper semi-analytic and all are $\R$--valued.
	
	Technically, \cite{BayZhou} uses the duality results in \cite{BN13} together with a minimax theorem to prove their results.
	Our first main result consists in introducing an enlarged canonical space (together with an enlarged canonical filtration)
	to reformulate the main problem as a superhedging problem for European options.
	Then by adapting the arguments in \cite{BN13}, we establish our duality under general conditions as in \cite{BN13}.
	Moreover, we do not assume that $\Phi_k$ is $\Fc_k$-measurable, which permits to study the superhedging problem for a portfolio containing an American option and some European options.
	Finally, our setting permits to use an approximation argument to study a new class of martingale optimal transport problem and to obtain a Kantorovich duality.

\section{A martingale (optimal) transport setup}
\label{sec:mot}

In this section we study the duality for American options in presence of a large family of statically traded European options. 
We assume that the statically traded options on the market are all vanilla options,
and are arbitrage--free (see \cite{CO} and \cite{CHO}) and numerous enough so that one can recover the marginal distribution of the underlying process $S$ at some maturity times $\T_0 = \{t_1, \cdots, t_{M}\} \subseteq \T$, where $t_{M} = N$.
More precisely, we are given a vector $\mu = (\mu_{1}, \cdots, \mu_{{M}})$ of marginal distributions. We write $\mu(\fu):=(\int \fu(x)\mu_1(dx),...,\int \fu(x)\mu_{M}(dx))$ and we assume that $\mu(|\cdot|) < \infty$ and
\be
\label{pcoc}
\mu_i(\fu) \le \mu_j(\fu)
		~~\mbox{for all}~
		i \le j, ~i,j \leq M,
		~\mbox{and convex function}~
		\fu: \R^d \to \R.
\ee
Here we work with $\Om :=\{s_0\}\times\R^{d\times N}$ where $s_0\in \R^d$, $S$ which is a canonical process on $\Om$ and $\Pc := \Bf(\Om)$.
Thus $\Omb := \Om \x \T$, $\Pcb=\frak P(\Omb)$.
The condition \eqref{pcoc} ensures the existence of a calibrated martingale measure, i.e. that the following sets are non--empty 
\begin{align*}
\Mc_\mu&:=\big\{\Q \in \Bf(\Om) ~:  \Law_{\Q}(S_{t_i})=\mu_i, ~i \leq M,~\mbox{and}~ S~\mbox{is an}~(\Q,\F) \mbox{-martingale} \big\},\\
\Mcb_\mu&:=	\big\{ 
			\Qb \in \Bf(\Omb) ~:  \Law_{\Qb}(S_{t_i})=\mu_i, ~i \leq M,~\mbox{and}~ S~\mbox{is a}~(\Qb,\Fbb) \mbox{-martingale} 
		\big\}.
\end{align*}
Let $\Lambda_0$ be the class of all Lipschitz functions $\lambda: \R^d \to \R$, 
and denote $\Lambda := \Lambda_0^{M}$. 
The statically traded options $g=(g^\no)_{\no\in\No}$ are given by $g^\no(\om):=\lambda(\om)-\mu(\lambda) $ where $\lambda(\om) := \sum_{i=1}^{M} \lambda_i(\om_{t_i})$ and 
$\mu(\lambda) := \sum_{i=1}^{M} \mu_{t_i}(\lambda_i)$.
Recall that $\Mc_g=\Mc_\mu$.
Since $\No$ is a linear space, the superhedging cost of the American option $\Phi$ using semi--static strategies $\pi^A_g(\Phi)$ defined in Subsection \ref{superhedging} can be rewritten as
\begin{align*}
\pi^A_g(\Phi)=\pi^A_\mu(\Phi)&:=\inf\{\mu(\lambda): \exists (H^1,...,H^N)\in \mathcal H^N \ \textrm{s.t.} \ H^j_i=H^k_i \ \forall 1\leq i\leq j\leq k\leq N\\ &\textrm{and} \ \lambda\in \Lambda\
\textrm{satisfying} \ \lambda(\om) \!+\! (H^k \!\is\! S)_N(\om) \geq \Phi_k(\om)
 \ \mbox{for all}~ k \in \T,~ \om \in \Om \}.
\end{align*}
Similarly, we denote by $\pib^E_{\mu}(\Phi)$ the corresponding superhedging cost for a European option with payoff $\Phi$ defined on $\Omb$, and one has $\pi^A_\mu(\Phi) = \pib^E_{\mu}(\Phi)$ by Theorem \ref{thm:American_European}.

\begin{Example} \label{exam:non_equiv}
Construct an example similar to Example \ref{ex:intro} to highlight that in \eqref{eq:PHAduality} we may have a strict inequality. Consider the case $N=2$, $\T_0 = \T = \{1,2\}$, $\mu_1 = \delta_{\{0\}}$ 
		and $\mu_2 = \frac{1}{4} \big( \delta_{\{-2\}} + \delta_{\{-1\}} + \delta_{\{1\}} + \delta_{\{2\}} \big)$.
		Let $\Phi_1(\{S_1=0\}) = 1$, $\Phi_2(\{|S_2|=1\}) = 2$ and $\Phi_2(\{|S_2|= 2\}) =0$.
		Then $\Mc_\mu$ contains only one probability measure $\Q$, and by direct computation, one has
		$$
			\E^{\Q} \big[ \Phi_{\tau}] = 1,
			~~~\mbox{for all}~~ \tau \in \Tc(\F).
		$$
		Let us now construct a martingale measure $\Qb_0$ by
		$$
			\Qb_0(d \om, d\theta) 
			~:=~
			\frac{1}{4} \delta_{\{1\}}(d \theta) \otimes \big( \delta_{(0,1)} + \delta_{(0,-1)} \big)(d \om)
			~+~
			\frac{1}{4} \delta_{\{2\}}(d \theta) \otimes \big( \delta_{(0,2)} + \delta_{(0,-2)} \big)(d \om).
		$$
		Then one can check that $\Qb_0 \in  \Mcb_\mu$ and it follows that
		$$
			\sup_{\Qb \in \Mcb_\mu}\E^{\Qb}[\Phi] \ge \E^{\Qb_0}[\Phi]
			~=~ \frac{3}{2}
			~>~ 1
			~=~ 
			\sup_{\Q \in \Mc_\mu} \sup_{\tau \in \Tc(\F)} \E^{\Q} \big[ \Phi_{\tau}].
		$$
Since the superhedging price of $\Phi$ equals to 3/2 as one can consider a superhedging strategy consisting of holding 3/2 in cash and one option $g$ from Example \ref{ex:intro}. 
In a similar way as in Example \ref{ex:intro} the duality may be recovered by allowing a dynamic trading options.
\end{Example}			

\subsection{Duality on the enlarged space $\Omb$}

The following theorem shows the duality for $\Omb$. Its proof is delegated to Section \ref{s:proofs3}.

\begin{Theorem} \label{theo:MOT}
Suppose that $\Phi : \Omb \to \R$ is bounded from above and upper semicontinuous.
Then there exists an optimal martingale measure $\Qb^* \in \Mcb_\mu$ and the pricing--hedging duality holds:
$$
\E^{\Qb^*}\big[\Phi\big]=\sup_{\Qb \in \Mcb_\mu}\E^{\Qb} \big[\Phi \big] ~~=~~ \pib^E_\mu(\Phi).
$$ 
and in particular \eqref{eq:PHAduality_Euro} holds. 
\end{Theorem}
\begin{Remark}
Note that in the above formulation each $\mu_i$ is an element of $\Bf(\R^d)$. Instead one could take $\mu_i$ to be an element of $(\Bf(\R))^d$. The same statements with analogous proofs would stay in force. This alternative formulation has more transparent financial interpretation since it corresponds only to marginal laws of terminal values of each stock price as opposed to the full distribution, see also \cite{Lim_mm} for a related discussion.
\end{Remark}

\subsection{Dynamic programming principle on $\Omh$}

\cite{Eldan} and
\cite{CoxKall} studied the Skorokhod embedding problem and the martingale optimal transport in continuous time using the measure--valued martingales. This point of view allows to obtain  the dynamic programming principle with marginal constraint since the terminal constraint is transformed into the initial constraint.
We adopt this perspective which proves to be very useful.\\
As before we work with the set of marginal times $\T_0=\{t_1, ..., t_{M}\}\subset\{1,...,N\}$ such that $t_{M}=N$, and marginal peacock measure $\mu=(\mu_{1},..., \mu_{{M}})$ where each $\mu_{i }$ is a probability measure on $\R^d$.
We let $\frak P_1(\R^d)=\{\eta\in \frak P(\R^d): \eta(|\cdot|)<\infty\}$ be the set of probability measures with finite first moment which we equip with the 1-Wasserstein distance, i.e. $\eta_n \to \eta_0$ if and only if
	$$
		 \int_{\R^d} \fu(x) \eta_n(dx) \to  \int_{\R^d} \fu(x) \eta_0(dx), ~~\forall \fu \in \Cc_1,	$$
where $\Cc_1$ denotes the set of all continuous functions on $\R^d$ with linear growth, which makes $\frak P_1(\R^d)$ a Polish space.
Continuing with the construction from Example \ref{ex:product}, $\Omh$ has to be an infinitely dimensional space and it is convenient to parametrize it as the canonical space for the measure--valued processes
$$\Omh := \{\mu\}\times(\frak P_1(\R^d))^{M\times {N}}$$ 
and denote  $\Xh=(\Xh^1_k, ..., \Xh^{M}_k)_{0\leq k\leq N}$ the canonical process on $\Omh$. Let $\Gh = (\Gch_k)_{0 \le k \le N}$ be the canonical filtration and $\Fh = (\Fch_k)_{0 \le k \le N}$ its universal completion. 
Denote by $\Tc(\Fh)$ the collection of all $\Fh$-stopping times.
For $\fu \in \Cc_1$ we denote the process of its integrals against $\Xh$ as
\begin{align*}
\Xh_k(\fu)&= (\Xh^1_k(\fu),...,\Xh^{M}_k(\fu)),\textrm{  where $\Xh^i_k(\fu):= \int_{\R^d} \fu(x) \Xh^{i}_k(dx)$ and}\\
\Xh_k(id) &=(\Xh^1_k(id),...,\Xh^{M}_k(id)), \textrm{  where $\Xh^i_k(id)= \int_{\R^d} x \Xh^{i}_k (dx)$}. 
\end{align*}
Define $\mathtt{i}: \Omh\to\Om$ by $\mathtt{i}(\omh)=(\Xh^{M}_0(id)(\omh),...,\Xh^{M}_N(id)(\omh))$ which is surjective and naturally extends processes on $\Om$ to processes on $\Omh$. 
In particular the price process extends via $S_k(\omh)=S_k(\mathtt{i}(\omh))=\Xh^{M}_k(id)(\omh)$ and the statically traded options via $g^\no(\omh)=g^\no(\mathtt{i}(\omh))=\no(\mathtt{i}(\omh))-\mu(\no)$.
Define a family of processes $\Y=(\Y^\no)_{\no\in\No}$ by $\Y^\no=\sum_{i=1}^{M}\Y^{\no_i}$ where 
\begin{align*}
\Y^{\no_i}_k&=
\begin{cases}
\Xh^i_k(\no_i)-\mu_i(\no_i)& 0\leq k\leq t_i-1\\
g^{\no_i}=\no_i(\Xh^i_{t_i}(id))-\mu_i(\no_i)&t_i\leq k\leq N
\end{cases}
\end{align*}
Note that $\Y^{\no_i}_0=0$.
\begin{Definition}
\rmi A probability measure $\Qh$ on $(\Omh, \Fch)$ is called a measure--valued martingale measure (MVM measure) if the process $(\Xh_k(\fu))_{0 \le k \le N}$ is a $(\Qh, \Fh)$-martingale for all $\fu \in \Cc_1$.\\
\rmii A MVM measure $\Qh$ is terminating if $\Xh^i_{t_i} \in \Delta:=\{\eta\in \frak P(\R^d): \eta=\delta_x, \ x\in \R^d\}$, $\Qh$-a.s.\\
\rmiii A MVM measure $\Qh$ is  consistent if $S_k=\Xh^i_k(id)$ for $k\leq t_i$ and $i=1, \cdots, M$, $\Qh$-a.s.
\end{Definition}
Let us denote by 
\begin{align*}
\Mch_\mu&=\{\Qh\in \Bf(\Omh\}: \Qh ~~\mbox{is terminating, consistent, MVM measure s.t.   } \forall \,i\leq M\; \}.
\end{align*}
The following lemma shows that the marginal distribution of $S$ at $t_i$ equals to $\mu_i$, $\Mch_\mu$-q.s. And hence $\Qh\circ \mathtt{i}^{-1}\in \Mc_\mu$ for any $\Qh\in\Mch_\mu$.
\begin{Lemma}
\label{cond:law}
For a measure $\Qh\in \Mch_\mu$ the following holds:\\
\rmi $\Law_{\Qh}(S_{t_i}| \Fch_k)=\Xh^i_k ~ \Qh\textrm{-a.s.}~\textrm{for} ~ k\leq t_i, ~ \textrm{and in particular}~ \Law_{\Qh}(S_{t_i})=\mu_i.$\\
\rmii For $k\leq t_j \leq t_i$, $\Xh^j_k\preceq \Xh^i_k$  $\Qh$-a.s., i.e., for any convex function $\fu$
$$\int_{\R^d}\fu(x)\Xh^j_k(dx)\leq \int_{\R^d}\fu(x)\Xh^i_k(dx) \quad \Qh\textrm{-a.s.}$$ 
\end{Lemma}

\proof
\rmi 
Let $A\subset \R^d$ and recall that $S_k=\Xh^i_k(id)$ $\Qh$-a.s.
Then we have
$$\int_{\R^d}\id_A(x)\Law_{\Qh}\left(\Xh^i_{t_i}(id) \Big| \Fch_k\right)(dx)
=\E^{\Qh}\left[\id_{\{\Xh^i_{t_i}(id)\in A\}} \Big|\Fch_k\right]
=\E^{\Qh}\left[\Xh^i_{t_i}(\id_A)\Big|\Fch_k\right]
=\Xh^i_k(\id_A),$$
where the second equality holds since $\Qh$ is terminating and the third one as $\Qh$ is MVM measure. Hence the first assertion is proven.\\
\rmii Let $j\leq i$, $k\leq t_j$ and $\fu$ be a convex function. 
Then
\begin{align*}
\int_{\R^d}\fu(x)\Xh^i_k(dx)&=
\E^{\Qh}\left[\fu\big(\Xh^i_{t_i}(id)\big)|\Fch_k\right]\\
&\geq \E^{\Qh}\left[\fu\big(\E^{\Qh}[\Xh^i_{t_i}(id)|\Fch_{t_j}]\big)|\Fch_k\right]\\
&=\E^{\Qh}\left[\fu\big(\Xh^i_{t_j}(id)\big)|\Fch_k\right]
=\E^{\Qh}\left[\fu\big(\Xh^j_{t_j}(id)\big)|\Fch_k\right]\\
&=\int_{\R^d}\fu(x)\Xh^j_k(dx)
\end{align*}
where the first and the last equality follow by \rmi, the penultimate is due to the consistency of $\Qh$ and the inequelity follows by conditional Jensen's inequality.
\qed

\vspace{2mm}

Let us denote by $\Mch$ the set of martingale measures for the dynamic extension $(\Omh, \Fch, \Fh, \mathtt{i}, Y)$ of $(\Om,\Fc,\F, \Pc,S,g)$ (see Section \ref{sec:dpp} for the definition of the dynamic extension). 

\begin{Lemma}
	\rmi Under any $\Qh\in\Mch_\mu$, the processes $S$ and $\Y^\no$, for $\no\in \No$, are $(\Qh, \Fh)$-martingales.
	In particular, one has $\Mch_\mu\subset\Mch$.
	
	\noindent \rmii The mapping $\mathtt{I}:\Mch_\mu\to\Mc_\mu$, 
	defined by $\mathtt{I}(\Qh)=\Qh\circ \mathtt{i}^{-1}$,
	is surjective. 
\end{Lemma}
	\proof \rmi The process $S=\Xh^{M}(id)$ is a $(\Qh, \Fh)$-martingale since $\Qh$ is MVM measure. 
To prove that $\Y^\no$ is a $(\Qh, \Fh)$-martingale for any $\no\in\No$, it is enough to show that for any $i\leq M$ and $\no \in \No_0$ one has $\E^{\Qh}\left[\no(\Xh^{M}_{t_i}(id))|\Fch_k\right]=\Xh^i_k(\no)$ for any $k<t_i$. 
The latter holds since
\begin{align*}
\E^{\Qh}\left[\no(\Xh^{M}_{t_i}(id))|\Fch_k\right]&=
\E^{\Qh}\left[\no(\Xh^{i}_{t_i}(id))|\Fch_k\right]=
\E^{\Qh}\left[\Xh^{i}_{t_i}(\no)|\Fch_k\right]=\Xh^i_k(\no),
\end{align*}
where the first equality follows by consistency of $\Qh$, the second since $\Qh$ is terminating and the last one as $\Qh$ is MVM measure. 

\vspace{1mm}

\rmii
Let $\Q\in \Mc_\mu$ and define the process $\eta=(\eta^1_k,.., \eta^{M}_k)_{k\leq N}$ by $\eta^i_k=\Law_{\Q}(S_{t_i}|\Fc_k)$. 
Note that by definition there exists a terminating, consistent  MVM measure $\Qh$ such that $\Qh[\Xh=\eta]=1$. 
\qed

\vspace{2mm}

For $\omh \in \Omh$, we define a set $[\omh]_{\Gch_k}$ as in \eqref{eq:om_k},
and denote by $\Mch^k_\mu(\omh)$ the following set of measures:
\begin{align*}
\Mch^k_\mu(\omh)&:=\Big\{\Qh\in \Bf(\Omh):\Qh~\mbox{is terminating and consistent,}\\
&~~~~~~~~~~~~~~~~~~~~~~~~~\Qh([\omh]_{\Gch_{k}})=1 ~~ \mbox{and $(\Xh_{l})_{k\leq l\leq N}$ is a $(\Qh,\Fh)$-MVM} \Big\}.
\end{align*}
Let us define a family of operators $\Ech_k$, etc., as in Example \ref{ex:non}:
\begin{align*}
\Ech_k(\xih)(\omh)&=\sup_{\Q\in\Mch^{k}_\mu(\omh)}\E^{\Q}\left[\xih\right],\quad \xih\in \fclh,
\end{align*}
and then the extension $\Ect_k$ as well as $\Ect^0$ on the enlarged space as in Section \ref{ameu}.
Then we have:
\begin{Theorem}
	For all upper semianalytic functionals $\xih : \Omh \to \R^N$,
	$\Ect_k(\xih)$ is also upper semianalytic and
\be \label{eq:claim_DPP_mvm}
\sup_{\Qt \in \Mct_{\mu}} \E^{\Qt} [ \xih] = \Ect^0 (\xih).
\ee
In particular the pricing--hedging duality \eqref{eq:PHA_dpp} holds in this MOT context for all functionals $\Phi: \Omb \to \R^N$ which are upper semicontinuous and bounded from above.
\end{Theorem}
\proof
	Notice that the pricing--hedging duality on $\Omb$ holds  by Theorem \ref{theo:MOT}.
	Then by Corollary \ref{cor:dpp_options},
	it is enough to establish the dynamic programming principle on $\Omh$ to prove the pricing-hedging duality \eqref{eq:PHA_dpp}.
	Using exactly the same arguments as in (4.12) of \cite{BN13}, to establish the dynamic programming principle on $\Omh$, it is enough to argue that $\Mch_\mu$ satisfies that
$$
	\{(\omh, \Qh) ~: \Qh \in \Mch_{\mu}^k(\omh) \}~~\mbox{is analytic}.
$$

	To prove the above analyticity property, we first observe that 
$$\Ech_k\circ ... \circ \Ech_{N-1}(\xih)(\omh)=\Ech_{k,k+1}\circ ... \circ \Ech_{N-1,N}(\xih)(\omh),\quad \xih\in \fclh$$
where $\Ec_{k,k+1}(\xih)(\om)=\sup_{\Mc^{k,k+1}_\mu(\om)}\E^{\Q}\left[\xih\right]$ and
\begin{align*}
\Mch^{k,k+1}_\mu(\omh)&:=\Big\{\Qh\in \Bf(\Omh):\Qh~\mbox{is terminating and consistent,}\\
&~~~~~~~~~~~~~~~~~~~~~~~~~~~~~\Qh([\omh]_{\Gch_{k}})=1 ~~ \omh_k(\fu) = \E^{\Qh}[ \Xh_{k+1}(\fu) ],
			~~\forall \fu \in \Cc_1\Big\}.
\end{align*} 	
	Next, let $\Cc_1^0$ denote a countable dense subset of $\Cc_1$ under the uniform convergence topology.
	Then it is clear that for each $k \in \T$, the set
\begin{align*} 
\big\{(\omh, \Qh)\in & \Omh \x \Bf(\Omh) :\  \Qh \in \Mch^{k,k+1}_\mu(\omh) \big\}
 =\big\{(\omh, \Qh) \in \Omh \x \Bf(\Omh) ~: \Qh([\omh]_{\Gch_{k}})=1\\
& ~~~~\Qh~\mbox{is terminating and consistent,} \ \ \omh_k(\fu) = \E^{\Qh}[ \Xh_{k+1}(\fu) ],
~~\forall \fu \in \Cc_1^0
\big\}
\end{align*}
	is a Borel set.
\qed

\section{Proofs for Section \ref{duality}}
\label{proofs1}

\proof[Proof of Proposition \ref{prop:dpp_Euro_American}]
First we prove that \eqref{eq:dpp_omb} implies \eqref{dpp_om}. For a given $\xi$ on $\Om$ let us define $\Psi$ on $\Omb$ by
\begin{align*}
\Psi((\om,k))&=-\infty\ \ \textrm{if} \ k\in \{1,...,N-1\}\\
\Psi((\om, N))&=\xi(\om) .
\end{align*}
Definition of $\Psi$ combined with \eqref{eq:dpp_omb} implies that 
\begin{align*}
\sup_{\Qb\in \Mcb}\E^{\Qb}[\Psi]&=\Ec_0\circ\Ec_1\circ...\circ\Ec_{N-1}\left(\xi\right).
\end{align*}
Moreover, one has that
\begin{align*}
\sup_{\Qb\in \Mcb}\E^{\Qb}[\Psi]&=\sup_{\Q\in\Mc}\E^{\Q}\left[\Psi_{N}\right]=\sup_{\Q\in\Mc}\E^{\Q}\left[\xi\right]
\end{align*}
since for a measure $\Qb\in \Mcb$ such that $\Qb(\Om\times \{1,...,N-1\})>0$ the expected value drops to $-\infty$.\\
Now let us prove that \eqref{dpp_om} and \eqref{inequality2} imply \eqref{eq:dpp_omb}. 
Define an $\F$-stopping time $\tau^*$ by 
\begin{align}
\label{tau}
\tau^*(\omega)&:=\inf \left\{k\geq 1: \Ec_k\left(\Psi(\cdot, k)\right)(\omega)= \Ecb^{k}(\Psi)(\om, k)\right\}\\
\nonumber
&= \inf \left\{k\geq 1: \Ec_k\left(\Psi(\cdot, k)\right)(\omega)\geq \Ec_k\left(\Ecb^{k+1}(\Psi)(\cdot, k+1)\right)(\omega)\right\}.
\end{align}
Note that on $\{k<\tau^*\}$  one has
\begin{align}
\label{before}
\Ec_k\left(\Psi(\cdot, k)\right)(\om)< \Ecb^{k}(\Psi)(\om, k)=\Ec_k\left(\Ecb^{k+1}(\Psi)(\cdot, k+1)\right)(\om).
\end{align}
Then 
\begin{align*}
\Ecb_0&\circ...\circ \Ecb_{N-1}(\Psi)=\\
&=\Ec_0\left(\id_{\{\tau^*=1\}}\Ec_1\left(\Ecb^{1}(\Psi)(\cdot, 1)\right)+\id_{\{\tau^*>1\}}\Ec_1\left(\Ecb^{1}(\Psi)(\cdot, 2)\right)\right)\\
&=...=\\
&=\Ec_0\circ\Ec_1\circ...\circ\Ec_{N-1}(\Psi_{\tau^*})\\
&=\sup_{\Q\in\Mc}\E^{\Q}\left[\Psi_{\tau^*}\right]
\end{align*}
where the last equality follows from DPP on $\Om$ \eqref{dpp_om}.
Note as well that 
\begin{align}
\label{ineaulity1}
\Ecb_0&\circ...\circ \Ecb_{N-1}(\Psi)=\sup_{\Q\in\Mc}\E^{\Q}\left[\Psi_{\tau^*}\right]\leq \sup_{\Q\in \Mc}\sup_{\tau\in \Tc(\F)}\E^{\Q}[\Psi_\tau]\leq
\sup_{\Qb\in \Mcb}\E^{\Qb}[\Psi].
\end{align}
Combining \eqref{inequality2}, we then conclude the proof.\qed

\begin{Proposition} \label{prop:non}
	The family $(\Ecb_k)$ given in Example \ref{ex:non} satisfies \eqref{inequality2}.
\end{Proposition}
\proof
In the context of Example \ref{ex:non}, the family $(\Ecb_k)$ take the following form:
\begin{align}
\Ecb_0(\Psi)&:=\sup_{\Q\in \Mc}\E^\Q[\Psi(\cdot, 1)]\\
\Ecb_k(\Psi)(\omb)&:=
\begin{cases}
\sup_{\Q\in \Mc^k(\om)}\E^\Q[\Psi(\cdot, \theta)]&\textrm{if }\theta<k\nonumber\\
\sup_{\Q\in \Mc^k(\om)}\E^\Q[\Psi(\cdot, k)]\lor \sup_{\Q\in \Mc^k(\om)}\E^\Q[\Psi(\cdot, k+1)]&\textrm{if } \theta\geq k.\nonumber
\end{cases}
\end{align}
To see that \eqref{inequality2} holds, it is insightful to rewrite $\Ecb^0$ in a slightly different way, as $\Ectt^0$ below. 
Let 
$$\Gcb^-_k:=\Gc_k\otimes \sigma(T\land k)\subset \Gcb_k:=\Gc_k\otimes \sigma(T\land (k+1))\subset \Fc_k\otimes \sigma(T\land (k+1))=:\Fcb_k,$$
$$\Mcb^{k,-}(\omb):=\{\Qb\lll\Pcb : \ \Qb\left[[\omb]_{\Gcb^-_k}\right]=1\ \textrm{and}\ \E^{\Qb}[\Delta S_n|\Fcb_{n-1}]=0 \ \forall n\in\{k+1,..., N\} \},$$
where $[\omb]_{\Gcb^-_k}$ is defined as in \eqref{eq:om_k}.
	Next, for $\Psi\in \fclb$, let us introduce operators 
\begin{align*}
\Ectt_0(\Psi):=\sup_{\Q\in \Mcb}\E^\Q[\Psi(\cdot, 1)],\qquad \Ectt_k(\Psi)(\omb):=
\sup_{\Qb\in \Mcb^{k,-}(\omb)}\E^{\Qb}[\Psi],\ k\leq N-1.
\end{align*}
Denote $ \Ecb^k(\cdot) := \Ecb_k \circ \cdots \circ \Ecb_{N-1}(\cdot)$ and $ \Ectt^k(\cdot) := \Ectt_k \circ \cdots \circ \Ectt_{N-1}(\cdot)$, and we claim that
\be \label{eq:dpp_two}
	\Ecb^k(\Psi)(\omb)=\Ectt^k(\Psi)(\omb),~0\leq k<N, ~\Psi\in \fclb.
\ee
Note that the conditional regular probabilities of any $\Q\in\Mcb$ w.r.t. $\Gcb^-_k$, denoted $\Q_{\omb}$, 
satisfy $\Q\left[\{\omb:\Q_{\omb}\in \Mcb^{k,-}(\omb)\}\right]=1$ and one has $\E^\Q[\Psi|\Fcb^-_k]\leq \Ectt_k(\Psi)$, $\Q$-a.s., which implies \eqref{inequality2} by the tower property of the conditional expectation and the definition of $\Ectt^0$.

Then it is enough to prove the claim \eqref{eq:dpp_two}.
Note that, for $\omb=(\om, \theta)$ with $\theta\leq k-1$, a measure $\Qb\in \Mcb^{k,-}(\omb)$ satisfies $\Qb_{\vert\Om}\in \Mc^k(\om)$ and $\Qb(\Om\times\{\theta\})=1$; and a measure $\Q\in \Mc_k(\om)$ satisfies $\Q\otimes \delta_\theta\in \Mcb^{k,-}(\omb)$. It is thus clear that, in this case, $\Ecb^k(f)(\omb)=\Ectt^k(f)(\omb).$

As a second step, for $\omb=(\om, \theta)$ with $\theta\geq k$, we show that $\Ecb_{k}(f)(\omb)\leq \Ectt_{k}(f)(\omb)$.
Take any $\Q\in\Mc^{k}(\om)$. Then, for $n\in\{k,..., N\}$, $\Q\otimes \delta_n\in \Mcb^{k,-}(\omb)$ and $\Q\otimes\delta_n(\Om\times \{n\})=1$. Hence it follows that $\Ecb_{k}(f)(\omb)\leq \Ectt_{k}(f)(\omb)$. 

In a final step, we show that, for $\omb=(\om, \theta)$ with $\theta\geq k$, $\Ecb_k(f)(\omb)\geq\Ectt_{k}(f)(\omb)$ holds. 
Let us start with $k=N-1$. 
Take any $\Qb\in\Mcb^{N-1,-}(\omb)$ and consider its r.c.p w.r.t. $\Gcb_{N-1}$ 
(the atom $\{\om\}\times\{N-1,N\}$ is divided into atoms $\{\om\}\times\{N-1\}$ and $\{\om\}\times\{N\}$) 
denoted by $\Qb_N$ and $\Qb_{N-1}$. Then, clearly, ${\Qb_N}_{\vert\Om}$ and ${\Qb_{N-1}}_{\vert\Om}$ belong to $\Mc^{N-1}(\om)$, 
and $\Qb_N(\{\om\}\times\{N\})=1$ and $\Qb_{N-1}(\{\om\}\times\{N-1\})=1$. Thus, it follows that $\Ecb_{N-1}(f)(\omb)\geq \Ectt_{N-1}(f)(\omb)$.

Finally, to complete the proof, we need to show that $\Ecb^{k+1}(f)(\omb)=\Ectt^{k+1}(f)(\omb)$ implies $\Ecb^k(f)(\omb)\geq\Ectt^k(f)(\omb)$ for $\omb=(\om, \theta)$ with $\theta\geq k$.
First note that $\Ecb^{k+1}(f)(\omb)=\Ectt^{k+1}(f)(\omb)$ is constant on $\theta\in \{k,...,N\}$, i.e., 
\begin{equation}
\label{constant}
\Ecb^{k+1}(f)((\om, \theta_1))=\Ecb^{k+1}(f)((\om,\theta_2))\ \ \ \textrm{for all $\om\in\Om$ and $\theta_1, \theta_2\in \{k,..., N\}$.}
\end{equation} 
Take any $\Qb\in\Mcb^{k,-}(\omb)$ and consider its r.c.p w.r.t. $\vartheta_N$ 
(the atom $\{\om\}\times\{k,...,N\}$ is divided into atoms $\{\om\}\times\{n\}$ for $n=k,...,N$) 
denoted by $\Qb_n$ for $n=k,...,N$. Then, clearly, ${\Qb_n}_{\vert\Om}\in\Mc^{k}(\om)$ 
and $\Qb_n([\om]_k\times\{n\})=1$ where $[\om]_k$ denotes an atom of $\Gc_k$ which contains $\om$. Thus, combining with \eqref{constant}, it follows that $\Ecb_{k}(f)(\omb)\geq \Ectt_{k}(f)(\omb)$.\qed

\section{Proofs for Section \ref{sec:main}}
\label{proofs2}

	We are now in the context of Section \ref{sec:main}, where $\Om_0 := \{\om_0\}$ is a singleton, $\Om_1$ is a nonempty Polish space and $\Om := \Om_0 \x \Om_1^N$.
	For technical reason, we introduce a $\Om_1$-valued canonical process $X = (X_k)_{0 \le k \le N}$ on the enlarged space $\Omb$ by
	$X_k(\omb) := \om_k$ for all $\omb := (\om, \theta) \in \Omb$, 
	and an enlarged filtration $\Gbb = (\Gcb_k)_{0 \le k \le N}$ by 
	$$
		\Gcb_0 := \{ \emptyset, \Omb\}
		~~~\mbox{and}~~
		\Gcb_k :=\sigma \big\{ X_i, ~\{T \le i\},~ i=1, \cdots, k \big\},
	$$
	and the universally completed filtration $\Fbb = (\Fcb_k)_{0 \le k \le N}$ by defining 
	$\Fcb_k$ as the universal completion of $\Gcb_k$.
	It follows that the random time $T: \Omb \to \T$ is an $\Gbb$-stopping time.
	We also define a restricted enlarged space, for every $k = 1, \cdots, N$,
	$$
		\Omb_k 
		~:=~
		\Om_k \x \{1, \cdots, k\} 
		~=~
		\Om_1^k \x \{1, \cdots, k\}.
	$$

\begin{Lemma}
	Let $\Pb\in\Pcb$ be a probability measure on $(\Omb, \Gcb_N)$,
	and $(\Pb_{\omb})_{\omb \in \Omb}$ be a family of regular conditional probability distribution of $\Pb$ w.r.t. $\Gcb_k$.
	Then for every $k =0, 1, \cdots, N-1$, one has $\Pb_{\omb} \circ X_{k+1}^{-1} \in \Pc_k(\om)$ for $\Pb$-a.e. $\omb = (\om, \theta) \in \Omb$.
\end{Lemma}

Let us introduce the following set of measures
\begin{eqnarray*}
		\Mcb_g^{loc} 
		\!\!\! &:=&\!\! 
		\{\Qb ~: \Qb \lll \Pcb, ~\E^{\Qb}[g^i] = 0,~i\in\{1,\cdots, e\} \\
		&&~~~~~~~~~~~~~~
		~\mbox{and}~
		S ~\mbox{is an}~ (\Fbb, \Qb) \mbox{-local martingale}
		\}.
\end{eqnarray*}

\begin{Lemma} \label{lemm:weak_duality}
Let $\Phi$ be upper semianalytic and $\Qb \in \Mcb^{loc}_g$.
Then for any $x \in \R$ and $(\Hb, h) \in \Hcb \x \R^e$ such that $x + (\Hb \is S)_N(\omb) + hg(\om) \ge \Phi(\omb)$, $\Qb$-a.s. one has
		$\E^{\Qb}\big[ \Phi \big] \le x$.
	\end{Lemma}
	\proof 
	The proof follows by exactly the same arguments as in Lemma A.2 of \cite{BN13}, using the discrete time local martingale characterization in Lemma A.1 of \cite{BN13}.
	\qed

\vspace{2mm}

Given $\Qb \in \Mcb^{loc}_0$ and $\varphi : \Om \to [0, \infty)$, we denote
\begin{eqnarray*}
\Mcb^{\varphi,\Qb} 
		\!\!\! &:=&\!\! 
		\{\Qb' \sim \Qb~: \E^{\Qb'}[\varphi] < \infty, ~\mbox{and}~
		S ~\mbox{is an}~ (\Fbb, \Qb') \mbox{-martingale}
		\}.
\end{eqnarray*}

Then by Lemma \ref{lemm:weak_duality}, one can easily obtain the weak duality:
\be\label{eq:weak_duality}
		\sup_{\Qb \in \Mcb_g} \E^{\Qb}[\Phi] 
		\le
		\sup_{\Qb \in \Mcb_g^{loc}} \E^{\Qb}[\Phi]
		\le
		{ \pib^E_g(\Phi)}.
\ee

\begin{Lemma} \label{lemm:red_Mcb_varphi}
Let $\Phi$ be upper semianalytic and $\Qb \in \Mcb^{loc}_0$ and $\varphi : \Om \to [1, \infty)$ be such that $|\Phi(\om,k)| \le \varphi(\om)$ for all $\omb=(\om,k) \in \Omb$.
		Then $\Mcb^{\varphi,\Qb} ~\neq~ \emptyset$,
		and moreover,
		$$\E^{\Qb}[\Phi]  ~\le~ \sup_{\Qb' \in \Mcb^{\varphi,\Qb}} \E^{\Qb'}[\Phi].$$
	\end{Lemma}
	\proof
		First, by Lemma 3.2 of \cite{BN13}, there exists a probability $\Pb_*$ equivalent to $\Qb$ on $(\Omb, \Fcb_N)$ such that $\E^{\Pb_*}[\varphi(X)] < \infty$.
		On the filtered probability space $(\Omb, \Fcb_N, \Fbb, \Pb_*)$, 
		one defines $\Mcb^{loc}_*$ as the collection of all probability measures $\Qb' \sim \Qb \sim \Pb_*$ under which $S$ is an $\Fbb$-local martingale.
Denote 
$$\pi^{E,\Qb}_0(\Phi): = \mathrm{inf} ~\big\{ x: \exists \Hb \in \Hcb ~~\mbox{s.t.}~~ x+ (\Hb \is S)_N \geq \Phi, \Qb\mbox{-a.s.} \big\} ,$$
		then by the classical arguments for the dominated discrete time market (such as  \cite{Kabanov, KabanovStricker}, see also Lemma A.3 of \cite{BN13}),
one can easily obtain the inequality 
$$
\E^{\Qb}[\Phi]	~\le~ \sup_{\Qb' \in \Mcb^{loc}_{*}} \E^{\Qb'}[\Phi] 
~\le~ \pib^{E,\Qb}_0(\Phi)
~\le~ 	\sup_{\Qb' \in \Mcb^{\varphi,\Qb}} \E^{\Qb'}[\Phi],
$$
which concludes the proof.
\qed

\vspace{2mm}

	Using Theorem 2.2 of \cite{BN13}, one can easily obtain a closedness result for the set of all payoffs which can be super--replicated from initial capital $x=0$, in our context.
	Let us denote by $\Lcb^0_+$ the set of all positive random variables on $\Omb$,
	and define 
	$$
		\Ccb ~:=~ \big\{ (\Hb \is S)_N + h g ~: \Hb \in \Hcb, h \in \R^e \big\} - \Lcb^0_+.
	$$
\begin{Lemma} \label{lemm:closedness}
Let $\Phi$ be upper semianalytic and $\NA(\Pcb)$. Then the set $\Ccb$ is closed in the following sense: \\
Let $(W^n)_{n \ge 1} \subset \Ccb$ and $W$ be a random variable such that $W^n \to W$, $\Pcb$-q.s., then $W \in \Ccb$.
\end{Lemma}
\proof
It is a direct consequence of Theorem 2.2 of \cite{BN13}, 
	where the results are given in a general abstract context.		
	\qed

\subsection{Proof of Theorem \ref{theo:main}: the case $e =0 $, equivalently $\No=\emptyset$}

For each $1 \le i \le j \le N$, we introduce an application from $\Om_j$ to $\Om_i$ (resp. $\Omb_j$ to $\Omb_i$) by
	$$
		[\om]_i \!:=\! (\om_1, \cdots, \om_i), ~\mbox{for all}~\om \in \Om_j
		~\mbox{(resp.}
		[\omb]_i ~:=~ ([\om]_i, \theta \wedge i),
		~\mbox{for all}~
		\omb = (\om, \theta) \in \Omb_j).
	$$	
Note that $\Fcb_k^-$ is the smallest $\sigma$-field on $\Omb$ generated by $[\cdot]_k: \Omb \to \Omb_k$;
	or equivalently, an $\Fcb_k^-$-measurable random variable $f$ defined on $\Omb$ can be identified as a Borel measurable function on $\Omb_k$.
	The canonical processes $X$ and $S$ are naturally defined on the restricted spaces 
	$\Om_k$ and $\Omb_k$.

	We next recall a notion of local no--arbitrage condition $\NA(\Pc_k(\om))$ introduced at the beginning of Section 4.2 in \cite{BN13}.
	Given a fixed $\om \in \Om_k$, we can consider $\Delta S_{k+1}(\om, \cdot) :=S_{k+1}(\om, \cdot) - S_k(\om)$ as a random variable on $\Om_1$, which determines a one--period market on $(\Om_1, \Bc(\Om_1))$ endowed with a class $\Pc_k(\om)$ of probability measures.
	Then $\NA(\Pc_k(\om))$ denotes the corresponding no--arbitrage condition in this one--period market, i.e., $\NA(\Pc_k(\om))$ holds if for all $H\in  \mathbb{R}^d$
		$$
			H \Delta S_{k+1}(\om, \cdot) \geq 0 
			\quad \Pc_k(\om)\mbox{-q.s.}~
			\Longrightarrow~ 
			H \Delta S_{k+1}(\om, \cdot) = 0
			\quad \Pc_k(\om)\mbox{-q.s.}
		$$

\begin{Lemma} \label{lemm:selectH}
		In the context of Section \ref{sec:main}, 
		let $f: \Omb_{k+1} \to \Rb$ be upper semianalytic,
		then $\Ecb_k(f): \Omb_k \to \Rb$ is still upper semianalytic.
		Moreover, there exist two universally measurable functions $(y_1, y_2): \Omb_k \to \R^d \x \R^d$
		such that
		\begin{eqnarray*}
			&&
			\Ecb_k(f) (\omb) ~+~ y_1(\omb) \Delta S_{k+1}(\om, \cdot)
			~\ge~ f(\om, \cdot, \theta)  
			~~\Pc_k(\om) \mbox{-q.s.}\\
			&&
			\Ecb_k(f) (\omb) ~+~ y_2(\omb) \Delta S_{k+1}(\om, \cdot)
			~\ge~ f(\om, \cdot, k+1) 
			~~\Pc_k(\om) \mbox{-q.s. }
		\end{eqnarray*}
		for all $\omb = (\om, \theta) \in \Omb_k$ such that $\NA(\Pc_k(\om))$ holds
		and $f(\om, \cdot, \theta) > - \infty, ~\Pc_k(\om)$-q.s.
		 $f(\om, \cdot, k+1) > - \infty, ~\Pc_k(\om)$-q.s.
	\end{Lemma}	
	
	\proof
		Notice that $f_1 \vee f_2$ is upper semianalytic whenever $f_1$ and $f_2$ are both upper semianalytic.
		Then the above lemma is an immediate consequence of Lemma 4.10 of \cite{BN13}
		as well as the definition of $\Ecb_k$.
	\qed
	
	\vspace{2mm}
	
Recall that $\Mcb_0$ (resp. $\Mcb^{loc}_0$) means $\Mcb_g$ (resp. $\Mcb^{loc}_g$) for the case $e=0$.

	\vspace{2mm}

	\noindent {\bf Proof of Theorem \ref{theo:main} (the case $e=0$).}
	First, one has the weak duality as in \eqref{eq:weak_duality}
	$$
		\sup_{\Qb \in \Mcb_g} \E^{\Qb}[ \Phi \big]
		~\le~
	\pib^E_g(\Phi).
	$$
Next, for the inverse inequality, we can assume, without loss of generality, that $\Phi$ is bounded from above.
	Indeed, by Lemma \ref{lemm:closedness}, one has
	$ \lim_{n \to \infty} \pib^E_g(\Phi \wedge n)  = \pib^E_g(\Phi)$
	(see also the proof of Theorem 3.4 of \cite{BN13}).
	Besides, the approximation 
	$\lim_{n \to \infty} \sup_{\Qb \in \Mcb_g} \E^{\Qb} [\Phi \wedge n] = \sup_{\Qb \in \Mcb_g} \E^{\Qb} [\Phi]$
	is an easy consequence of the monotone convergence theorem.
	
	When $\Phi$ is bounded from  above, by Lemma \ref{lemma:super_hedging}, it is enough to prove that there is some $\Hb \in \Hcb$ such that
	\be \label{eq:superheding_Ec}
		\Ecb^0[\Phi]
		~+~
		(\Hb \is S)_N
		~\ge~
		\Phi
		~~~\Pcb\mbox{-q.s.}
	\ee
In view of Lemma \ref{lemm:red_Mcb_varphi}, we know $\Ecb^k(\Phi)(\omb) > -\infty$ for all $\omb \in \Omb_k$.
	Further, by Lemma \ref{lemm:selectH},
	there exist two universally measurable functions 
	$(y^k_1, y^k_2): \Omb_k \to \R^d \x \R^d$ such that
	\begin{eqnarray*}
		&&
		y^k_1(\omb) \Delta S_{k+1}(\om, \cdot) \ge \Ecb^{k+1}(\Phi)(\om, \cdot, \theta) - \Ecb^k(\Phi)(\omb)
		~~\Pc_k(\om) \mbox{-q.s.} \\
		&&
		y^k_2(\omb) \Delta S_{k+1}(\om, \cdot) \ge \Ecb^{k+1}(\Phi)(\om, \cdot, k+1) - \Ecb^k(\Phi)(\omb)
		~~\Pc_k(\om) \mbox{-q.s.}
	\end{eqnarray*}
	for all $\omb = (\om, \theta) \in \Omb_k$ such that $\NA(\Pc_k(\om))$ holds.

	Since $N_k := \{\om_k : \NA(\Pc_k(\om)) ~\mbox{fails}\}$ is $\Pc$-polar by Theorem 4.5 of \cite{BN13},
	it follows that, with $\Hb_{k+1}(\omb) := y^k_1([\omb]_k) \1_{\{\theta \le k\}} + y^k_1([\omb]_k) \1_{\{\theta > k\}}$, one has
	$$
		\sum_{k=0}^{N-1} \Hb_{k+1} \Delta S_{k+1} 
		~\ge~ 
		\sum_{k=0}^{N-1} \big(\Ecb^{k+1}(\Phi) - \Ecb^k(\Phi) \big)
		~=~
		\Phi - \Ecb(\Phi),
		~~~\Pc \mbox{-q.s.}
	$$
	To conclude, it is enough to notice that the above $\Hb$ is an optimal dual strategy for the case $\Phi$ being bounded from above.
	The existence of the optimal dual strategy for general $\Phi$ is then a consequence of Lemma \ref{lemm:closedness}.
	\qed

\subsection{Proof of Theorem \ref{theo:main}: the case $e \ge 1$, equivalently $\No\neq\emptyset$}
	
	We will adapt the arguments in Section 5 of \cite{BN13} to prove Theorem \ref{theo:main} in the context with finitely many options $e \ge 1$.
	
	For technical reasons, we introduce 
	$$ 
		\varphi(\om, \theta) := 1 + |g^1(\om)| + \cdots + |g^e(\om)| + \max_{1 \le k \le N} |\Phi_k(\om)|,
	$$ 
	which depends only on $\om$, and
	\be
		\Mcb^{\varphi}_g
		:=
		 \{ 
		 	\Qb \in \Mcb_0 ~: \E^{\Qb}[\varphi] < \infty ~\mbox{and}~
		 	\E^{\Qb}[g^i] = 0 
		 	~\mbox{for}~ i=1, \cdots, e 
		 \}.
	\ee
	Moreover, in view of Lemma \ref{lemm:red_Mcb_varphi}, one has
	$$\sup_{\Qb \in \Mcb_g} \E^{\Qb}[ \Phi] = \sup_{\Qb \in \Mcb_g^{\varphi}} \E^{\Qb}[ \Phi].$$

	\noindent {\bf Proof of Theorem \ref{theo:main} (the case $e \ge 1$)}.
	The existence of some $\Qb \in \Mcb_g$ is an easy consequence of Theorem 5.1 of \cite{BN13} under $\NA(\Pc)$.
	Moreover, similarly to \cite{BN13}, there exists an optimal dual strategies by Lemma \ref{lemm:closedness}.
	We will then focus on the duality results.
	
First, the duality ${\pib^E_g(\Phi)} = \sup_{\Qb \in \Mcb_g} \E^{\Qb}[\Phi]$ in \eqref{eq:main} has already been proved for the case $e=0$, 
	we will use the induction arguments:
	Suppose that the duality \eqref{eq:main} holds true for the case with $e \geq 0$,
	We aim to prove the duality with $e+1$ options:
	$$
		\pib^E_{(g,f)} (\Phi) 
		~=~ 
		\sup_{\Qb \in \Mcb_{(g,f)}^{\varphi}} \E^{\Qb}[\Phi],
	$$	
	where the additional option has a Borel--measurable payoff function $f \equiv g^{e+1}$ 
	such that $|f| \leq \varphi$, and has an initial price $f_0 = 0$.
	By the weak duality in \eqref{eq:weak_duality} and Lemma \ref{lemm:red_Mcb_varphi}, the ``$\ge$'' side of the inequality holds true, 
	we will focus on the ``$\le$'' side of the inequality:
	\be \label{eq:ineq_e_1}
		\pib^E_{(g,f)} (\Phi) 
		~\le~ 
		\sup_{\Qb \in \Mcb_{(f,g)}^{\varphi}} \E^{\Qb}[\Phi].
	\ee

	If $f$ is replicable by some semi--static strategy with underlying $S$ and options $(g^1, \cdots, g^e)$ in sense that $\exists \Hb \in \Hcb, h \in \R^e$, s.t. $f=(\Hb \is S)_N + hg, \Pcb\mbox{-q.s.}$
	(or equivalently, $\exists H \in \Hc, h \in \R^e$, s.t. $f=(H \is S)_N + hg, \Pc\mbox{-q.s.}$), 
	then the problem is reduced to the case with $e$ options and the result is trivial. 
	Let us assume that $f$ is not replicable, and we claim that there exists a sequence $(\Qb_n)_{n \ge 1} \subset \Mcb_g^{\varphi}$ such that
	\be \label{eq:claim_e_1}
		\E^{\Qb_n} [f] ~\longrightarrow~ f_0
		~~~\mbox{and}~~ 
		\E^{\Qb_n} [\Phi] ~\longrightarrow~ \pib^E_{(g,f)} (\Phi),
		~~~\mbox{as}~~
		n \longrightarrow \infty.
	\ee
Next, denote by $\pib^E_g(f)$ the minimum superhedging cost of European option $f$ using $S$ and $(g^1, \cdots, g^e)$, i.e. 
$$\pi^E_g(f)=\pib^E_g(f) = \mbox{inf}~ \{ x: \exists \Hb \in \Hcb, h \in \R^e, \mbox{s.t.}~ x+(\Hb \is S)_N + hg \ge f, ~~\Pcb\mbox{-q.s.} \}.$$
	Since $f$ is not replicable, by the second fundamental theorem in Theorem 5.1.(c) of \cite{BN13}, we have that $\Qb \mapsto \E^{\Qb}[f]$ is not constant on $\Mcb_g^{\varphi}$. 
	Then, under the no--arbitrage condition, one has $0 = f_0 <  \pib^E_g(f)$.
	It follows that
	$0 = f_0 <  \pib^E_g(f) = \sup_{\Qb \in \Mcb_g^{\varphi}} \E^{\Qb} [f(]$.
	Thus there exists some $\Qb_+ \in \Mcb_g^{\varphi}$, s.t. $0 < \E^{\Qb_+}[f] < \pib^E_g(f)$.
With the same argument on $-f$, we can find another $\Qb_- \in \Mcb_g^{\varphi}$ such that
	$$
		-\pib^E_g(-f) < \E^{\Qb_-} [f] < f_0 < \E^{\Qb_+} [f] < \pib^E_g(f)
	$$
	Then one can choose an appropriate sequence of weight 
	$(\lambda_-^n, \lambda_0^n, \lambda_+^n) \in \R_+^3$, such that $\lambda_-^n + \lambda^n_0 + \lambda_+^n = 1$, 
	$\lambda_{\pm}^n \rightarrow 0$
	and
	$$
		\Qb_n':= \lambda_-^n \Qb_- + \lambda^n_0 \Qb_n + \lambda_+^n \Qb_+ \in \Mcb_g,
		~~\mbox{and}~
		\E^{\Qb_n'} [f] = f_0 = 0,
	$$
	i.e. $\Qb_n' \in \Mcb_{(g,f)}^{\varphi}$.
	Moreover, since $\lambda_{\pm}^n \rightarrow 0$, it follows that
	$\E^{\Qb_n'}[\Phi] \to \pib^E_{(g,f)}(\Phi)$ 
	and we hence have the inequality \eqref{eq:ineq_e_1}.
	
	It is enough to prove the claim \eqref{eq:claim_e_1},
	for which we suppose without loss of generality that $\pib^E_{(g,f)}(\Phi) = 0$. 
	Assume that \eqref{eq:claim_e_1} fails, then one has
	$$
		    0 ~\notin~ \overline{ \{ \E^{\Qb} [ (f,\Phi) ]~: \Qb \in \Mcb_g^{\varphi}  \} } 
		    ~\subseteq~
		    \R^2.
	$$
By the convexity of the above set and the separation argument, 
	there exists $(y,z) \in \mathbb{R}^2$ with $|(y,z)|=1$, such that
\be\label{eq:contradition_e_1}
0~>~\sup_{\Qb \in \Mcb_g^{\varphi}} \E^{\Qb} [ yf+z \Phi]~=~\pib^E_g(yf + z \Phi)~\ge~\pib^E_{(g,f)}(z \Phi).
\ee
The strict inequality $\pib^E_{(g,f)}(z \Phi) < 0$ implies that $z \neq 0$.
Now, if $z>0$, we then have $\pib^E_{(g,f)} (\Phi) < 0$, which contradicts $\pib^E_{(f,g)} (\Phi) =0$. 
	If $z<0$, then by \eqref{eq:contradition_e_1}, one has $0 > \E^{\Qb'} [ yf+z \Phi] = \E^{\Qb'} [ z \Phi ]$ for some $\Qb' \in \Mcb_{(g,f)} \subseteq \Mcb_g$ since $\Mcb_{(f,g)}$ is nonempty under the $\NA(\Pcb)$ assumption in the case of $e+1$ options.
Then in the case $z<0$, one has $\E^{\Qb'} [ \Phi ] > 0 = \pib^E_{(f,g)}(\Phi)$, 
	which contradicts the weak duality result \eqref{eq:weak_duality},
	and we hence conclude the proof of the duality.
	\qed

\section{Proofs for Section \ref{sec:mot}}
\label{s:proofs3}

A first idea how to prove Theorem \ref{theo:MOT} could be the following two steps argument as in  \cite{GTT2}.
Firstly, under the condition that $\Phi$ is bounded from above and upper semicontinuous, one could prove that 
$$\mu \in \Bf((\R^d)^{M}) \mapsto \sup_{\Qb \in \Mcb_\mu}\E^{\Qb} \big[\Phi \big] \in \R$$ 
is concave and upper semicontinuous, 
where we equip $\Bf((\R^d)^{M})$ with a Wasserstein kind topology.
Secondly, using Fenchel--Moreau theorem, it follows that
\be \label{eq:duality_PD0}
\sup_{\Qb \in \Mcb_\mu}\E^{\Qb} \big[\Phi \big] 
		~=~
		\pib^E_{\mu,0}(\Phi)
		~:=~
		\inf_{\lambda \in \Lambda} \Big\{\mu(\lambda) 
			+ \sup_{\Qb \in \Mcb_0} \E^{\Qb} \big[ \Phi - \lambda \big] 
		\Big\}.
	\ee
Solving the maximization problem \eqref{eq:duality_PD0}, by using Theorem \ref{theo:main}, concludes the proof of Theorem  \ref{theo:MOT}.

	However in the following, we will provide another proof, which is based on an approximation argument.
	For simplicity, we suppose that $\T_0 = \{N\}$, where the same arguments work for more general $\T_0$.
	In preparation, let us provide a technical lemma.
	In the context  of the martingale optimal transport problem,
	we introduce a sequence of payoff functions $(g^i)_{i \ge 1}$ by
	$$
		g^i(\om) := f^i(\om_N) - c^i
		~~~\mbox{with}~~
		c^i:= \int_{\R^d} f^i(x) \mu(dx),
	$$
	where $f^i : \R^d \to \R$ is Lipschitz and $(f^i)_{i \ge 1}$ is dense in the space of all Lipschitz functions on $\R^d$ under the uniform convergence topology,
	and moreover, it contains all functions in form $(x_j -K_n)^+$, $(-K_n - x_j)^+$ for $j=1, \cdots, d$ and $n \ge 1$, where $(K_n)_{n \ge 1} \subset \R$ is a sequence such that $K_n \to \infty$.
	Notice that $\mu$ has finite first order and hence $c^i$ are all finite constants.

Next, let us introduce an approximate dual problem by
\begin{eqnarray*}
		\pi^A_{\mu,m}(\Phi)
		&:=&
		\inf \Big\{
			x ~: \exists(\Hb, h) \in \Hcb \x \R^m ~\mbox{s.t.}~ 
			\mbox{for all}~ k \in \T, ~\om \in \Om,\\
		&&~~~~~~~~~~~~
			x + \sum_{i=1}^m h^i g^i(\om_N) + (\Hb^k \is S)_N(\om) \ge \Phi_k(\om)
		\Big\}.
\end{eqnarray*}

Similarly,
	$$
		\Mcb_{\mu,m} ~:=~
		\big\{ 
			\Qb \in \Mcb ~: \E^{\Qb}[g^i] = 0 ~~\mbox{for}~i = 1, \cdots, m
		\big\},
	$$
	and
	$$
		P_{\mu,m} ~:=~ \sup_{\Qb \in \Mcb_{\mu,m}} \E^{\Qb} \big[ \Phi \big].
	$$

	\begin{Lemma}\label{lemm:compact1}
		Let $(\Qb_m)_{m \geq 1} \subset \Mcb$ be a sequence of martingale measures 
		such that $\Qb_m \in \Mcb_{\mu,m}$ for each $m \ge 1$.
		Then,

		\noindent \rmi $(\Qb_m)_{m \geq 1}$ is relatively compact under the weak convergence topology.
		
		\noindent \rmii The sequence $(S_N, \Qb_m)_{m \ge 1}$ is uniformly integrable,
		and any accumulation point of $(\Qb_m)_{m \geq 1}$ belongs to $\Mc_\mu$.
	\end{Lemma}
	\proof \rmi Without loss of generality, we assume that $f_1(x) = \sum_{i=1}^d |x_i|$ so that	
	$$\sup_{m \ge 1} \E^{\Qb_m} \big[ \sum_{i=1}^d |S_N^i| \big] < \int_{\R^d} \sum_{i=1}^d|x_i| \mu(dx) < \infty.$$
	Let us first prove the relative compactness of $(\Qb_m)_{m \geq 1}$.
	By Prokhorov theorem, it is enough to find, for every $\eps > 0$, 
	a compact set $D_{\eps} \subset \R^d$ such that $\Qb_m[ S_k \notin D_{\eps} ] \le \eps$ for all $k = 1, \cdots, N$.
	It is then enough to find, for every $\eps > 0$, a constant $K_{\eps} > 0$ such that
	$\Qb_m \big[ |S^i_k| \ge K_{\eps} \big] \le \eps$ for all $i = 1, \cdots, d$ and $k=1, \cdots, N$.
	Next, by the martingale property, one has $\E^{\Qb_m} [|S_k^i|] \le \E^{\Qb_m}[|S_N^i|]$.
	Then for every $\eps>0$, one can choose $K_{\eps} > 0$ such that $\sup_{m \ge 1} \E^{\Qb_m} \big[ \sum_{i=1}^d |S_N^i| \big] \le K_{\eps} \eps.$
	It follows that $\Qb_m \big[ |S^i_k| \ge K_{\eps} \big] \le \frac{\E^{\Qb_m}[|S^i_k|]}{K_{\eps}} \le  \eps$,
	and hence $(\Qb_m)_{m \ge 1}$ is relatively compact.
	\vspace{1mm}
	
	\noindent \rmii To see that the sequence $(S_N, \Qb_m)_{m \ge 1}$ is uniformly integrable,
	it is enough to notice that $|x_i| \1_{|x_i| \ge 2K_n} \le 2 (|x_i| - K_n) \1_{|x_i| \ge K_n}$, where the latter is a payoff function contained in the sequence $(f_k)_{k \ge 1}$.
	
	\vspace{1mm}

	\noindent \rmiii Let $\Qb_0$ be an accumulation point of  $(\Qb_m)_{m \geq 1}$.
	Since the sequence $(f_k)_{k \ge 1}$ is supposed to be dense in the space of all Lipschtiz functions on $\R^d$ under the uniformly convergence topology,
	it is easy to obtain that $\Qb_0 \circ S_N^{-1} = \mu$.

	\vspace{1mm}

	\noindent \rmiv To conclude the proof,
	it is enough to show that the martingale property is preserved for the limiting measure $\Qb_0$.
	By abstracting a subsequence, we assume that $\Qb_m \to \Qb_0$ weakly,
	and we will prove that for all $1 \leq k_1 < k_2 \leq N$, 
	for any bounded continuous function $\varphi: (\R^d)^{k_1} \x \T \rightarrow \R$,
	one has
	\begin{equation} \label{eq:martinglim}
		\E^{\Qb_0} \big[ \varphi \big(S_1, \cdots, S_{k_1}, T \wedge (k_1+1) \big) (S_{k_2}-S_{k_1}) \big]
		~=~
		0.
	\end{equation}
	Let $K >0$, and $\chi_K: \R^d \to \R^d$ a continuous function uniformly bounded by $K$ satisfying
	$\chi_K(x) = x$ when $\|x\| \le K$, and $\chi_K(x) = 0$ when $\|x \| \ge K+1$.
	Then for every $m =0$ or $m \ge 1$, one has
	\be \label{eq:mart_ineq_1}
		\big| \E^{\Qb_m} \big[ \varphi(S,T) (S_{k_2}-S_{k_1}) \big] \big|
		\!\!&\le &\!\!
		\big| \E^{\Qb_m} \big[ \varphi(S,T) \big(\chi_K(S_{k_2}) - \chi_K(S_{k_1}) \big) \big] \big| \nonumber\\
		\!\!&&\!\!\!\!\!\!+
		|\varphi|_{\infty} 
		\E^{\Qb_m} \big[ |S_{k_2}|\mathbf{1}_{|S_{k_2}| \geq K}+|S_{k_1}|\mathbf{1}_{|S_{k_1}| \geq K} \big],
		~~~~~
	\ee
	where we simplify $\varphi(S_1, \cdots, S_{k_1}, T \wedge (k_1+1))$ to $\varphi(S, T)$.

	For every $\eps > 0$, by uniformly integrability of $(S_N, \Qb_m)_{m \ge 1}$, there is $K_{\eps} > 0$
	such that 
	\be \label{eq:mart_ineq_2}
		|\varphi|_{\infty} 
		\E^{\Qb_m} 
		\big[ |S_{k_2}|\mathbf{1}_{|S_{k_2}| \geq K_{\eps}}
			+|S_{k_1}|\mathbf{1}_{|S_{k_1}| \geq K_{\eps}} 
		\big]
		~\le~
		\eps,
		~~~\mbox{for all}~m=0, 1, \cdots
	\ee
	Moreover, for $m \ge 1$, $\Qb_m$ is a martingale measure,
	then $\E^{\Qb_m} \big[ \varphi(S,T) (S_{k_2}-S_{k_1}) \big] = 0$ and hence
	$\big| \E^{\Qb_m} \big[ \varphi(S,T) \big(\chi_K(S_{k_2}) - \chi_K(S_{k_1}) \big) \big] \big| \le \eps$.
	Then by taking the limit $m \rightarrow \infty$, it follows that
	\be \label{eq:mart_ineq_3}
		\big| \E^{\Qb_0} \big[ \varphi(S,T) \big(\chi_K(S_{k_2}) - \chi_K(S_{k_1}) \big) \big] \big| 
		~\le~
		\eps.
	\ee
	Combining \eqref{eq:mart_ineq_1}, \eqref{eq:mart_ineq_2} and \eqref{eq:mart_ineq_3}, and by the arbitrariness of $\eps > 0$,
	it follows that \eqref{eq:martinglim} holds true and we hence conclude the proof.
	\qed

	\vspace{3mm}

	\noindent {\bf Proof of Theorem \ref{theo:MOT}}.
	We notice that by Theorem \ref{theo:main}, 
$$\sup_{\Qb \in \Mcb_{\mu,m}} \E^{\Qb} \big[ \Phi \big] = \pi^A_{\mu,m}(\Phi) \ge \pi^A_\mu(\Phi).$$
	Let $(\Qb_m)_{m \ge 1}$ be a sequence of probability measures such that $\Qb_m \in \Mcb_{\mu,m}$ for each $m \ge 1$ and
	$$
		\limsup_{m \to \infty} \E^{\Qb_m} \big[ \Phi_T(S) \big] 
		~=~
		\limsup_{m \to \infty} \sup_{\Qb \in \Mcb_{\mu,m}} \E^{\Qb} \big[ \Phi \big].
	$$
	It follows by Lemma \ref{lemm:compact1} that
	there is some $\Qb_0 \in \Mcb_\mu$ and a subsequence $\Qb_{m_k} \to \Qb_0$ 
	under the weak convergence topology. 
	Using upper semi--continuity of $\Phi$ and by Fatou's lemma, it follows that $\E^{\Qb_0} \big[ \Phi_T(S) \big] \ge \limsup_{m \to \infty} \E^{\Qb_m} \big[ \Phi_T(S) \big]$.
	It leads to the inequality
	$$
		 \sup_{\Qb \in \Mcb_\mu}\E^{\Qb} \big[\Phi \big] 
		~\ge~ 
		\E^{\Qb_0} \big[ \Phi \big] 
		~\ge~
		\limsup_{m \to \infty} \sup_{\Qb \in \Mcb_{\mu,m}} \E^{\Qb} \big[ \Phi \big]
		 ~=~
		\limsup_{m \to \infty} \pi^A_{\mu,m}(\Phi)
		~\ge~
		\pi^A_\mu(\Phi)	
	$$
	and we hence conclude the proof by the weak duality \eqref{eq:weak_duality}.
	\qed

\bibliographystyle{apalike}
\bibliography{biblio}

\end{document}